\documentclass[onecolumn]{IEEEtran}
\usepackage{cite}

\usepackage{graphicx,epsfig,psfrag}
\usepackage{epstopdf}
\usepackage{graphicx,amsmath,amssymb}
\usepackage{overpic}
\usepackage{amsfonts}
\usepackage{graphicx}
\usepackage{subfigure}
\usepackage{ifthen}
\usepackage{color}
\usepackage{fancyhdr}
\usepackage{stfloats}

\usepackage{multirow}
\usepackage{algorithm}
\usepackage{algorithmic}

\usepackage{picins}

\usepackage{threeparttable}

\newcommand{\R}{{\rm I\!R}}

\newcommand{\diag}{\operatorname{diag}} % diagonal part
\usepackage{color}

\newcommand{\tabincell}[2]{\begin{tabular}{@{}#1@{}}#2\end{tabular}}
\newcommand{\mytoprule}{\rule{\linewidth}{1pt}\vspace{-1ex}}

\usepackage{caption}
\def\bb#1{{\mathbb #1}}

\newtheorem{theorem}{Theorem}
\newtheorem{proposition}{Proposition}
\newtheorem{corollary}{Corollary}
\newtheorem{definition}{Definition}
\newtheorem{lemma}{Lemma}
\newtheorem{remark}{Remark}
\newtheorem{assumption}{Assumption}
\newtheorem{exmp}{Example}

\begin{document}
%
% paper title
% can use linebreaks \\ within to get better formatting as desired
\title{Optimizing Pinning Control of Complex Dynamical Networks Based on Spectral Properties of Grounded Laplacian Matrices}

\author{Hui~Liu,~\IEEEmembership{Member,~IEEE,}
        Xuanhong Xu,~  % <-this % stops a space
        Jun-An~Lu,~  % <-this % stops a space
        Guanrong Chen,~\IEEEmembership{Fellow,~IEEE,}
        \\and Zhigang Zeng,~\IEEEmembership{Senior~Member,~IEEE}
%\thanks{Manuscript received April 19, 2005; revised January 11, 2007.}
%\thanks{This work was submitted to the journal of IEEE Transactions on Systems, Man, and Cybernetics: Systems.
%}
%\thanks{}% <-this % stops a space
\thanks{H. Liu, X. Xu, and Z. Zeng are with the School of Automation
and also with the Key Laboratory of Image Processing and Intelligent Control of Education Ministry of China, Huazhong University of Science and Technology, Wuhan 430074, China.  {\tt\small hliu@hust.edu.cn; xhxu@hust.edu.cn; zgzeng@hust.edu.cn}}
\thanks{J.-A.~Lu is with the School of Mathematics and Statistics, Wuhan University, Wuhan 430072, China. {\tt\small jalu@whu.edu.cn}}
\thanks{G. Chen is with the Department of Electronic Engineering, City University of Hong Kong, Hong Kong Special Administrative Region, China  {\tt\small eegchen@cityu.edu.hk}}
%\thanks{}% <-this % stops a space]
}
%\thanks{Manuscript received April 19, 2005; revised January 11, 2007.}

%\thanks{}% <-this % stops a space

% note the % following the last \IEEEmembership and also \thanks -
% these prevent an unwanted space from occurring between the last author name
% and the end of the author line. i.e., if you had this:
%
% \author{....lastname \thanks{...} \thanks{...} }
%                     ^------------^------------^----Do not want these spaces!
% make the title area
\maketitle

\begin{abstract}
Pinning control of a complex network aims at forcing the states of all nodes to track an external signal by controlling a small number of nodes in the network. In this paper, an algebraic graph-theoretic condition is introduced to optimize pinning control. When individual node dynamics and coupling strength of the network are given, the effectiveness of pinning scheme can be measured by the smallest eigenvalue of the grounded Laplacian matrix obtained by deleting the rows and columns corresponding to the pinned nodes from the Laplacian matrix of the network. The larger this smallest eigenvalue, the more effective the pinning scheme. Spectral properties of the smallest eigenvalue are analyzed using the network topology information, including the spectrum of the network Laplacian matrix, the minimal degree of uncontrolled nodes, the number of edges between the controlled node set and the uncontrolled node set, etc. The identified properties are shown effective for optimizing the pinning control strategy, as demonstrated by illustrative examples. Finally, for both scale-free and small-world networks, in order to maximize their corresponding smallest eigenvalues, it is better to pin the nodes with large degrees when the percentage of pinned nodes is relatively small, while it is better to pin nodes with small degrees when the percentage is relatively large. This surprising phenomenon can be explained by one of the theorems established.
\end{abstract}
% IEEEtran.cls defaults to using nonbold math in the Abstract.
% This preserves the distinction between vectors and scalars. However,
% if the journal you are submitting to favors bold math in the abstract,
% then you can use LaTeX's standard command \boldmath at the very start
% of the abstract to achieve this. Many IEEE journals frown on math
% in the abstract anyway.

% Note that keywords are not normally used for peerreview papers.

% For peer review papers, you can put extra information on the cover
% page as needed:
% \ifCLASSOPTIONpeerreview
% \begin{center} \bfseries EDICS Category: 3-BBND \end{center}
% \fi
%
% For peerreview papers, this IEEEtran command inserts a page break and
% creates the second title. It will be ignored for other modes.
\IEEEpeerreviewmaketitle

\section{INTRODUCTION}
Controlling a complex network to achieve a certain desired objective is an important task for various interacting systems, regarding such as trajectory tracking in a group of mobile autonomous agents \cite{ChLiZhChXi16,XiYuChGa16}, generators synchronization in power grids \cite{DoChBu13}, and rumor propagation in social networks \cite{Za02}. In practical situations, it is expensive and unpractical to control all nodes especially in a large-scale complex network.

To overcome the difficulties, the concept of pinning control was proposed in \cite{WaCh02a,LiWaCh04}. It aims at forcing the states of all nodes to a desired trajectory by controlling only a small portion of nodes in a given network.
From a theoretical viewpoint, there have been significant achievements in the research of pinning control.
By considering a specific trajectory as the state of a virtual dynamical system, pinning control can also be understood as a synchronization problem.
A number of conditions and criteria for pinning synchronization have been proposed in the last decade \cite{WaCh02a,LiWaCh04,ZhLuLu08,YuChLu09,Wu08a,Wu08b,SodiGaCh07,WaDaDoCaSu08,SoCa10,ChLilu07,Ch14,YuChLuKu13,SuRoChWaChWa13,
WaWuHuReWu16a,ChYuSu17,YuChWaYa09,XiHoWa11,LiCh15a}, to name just a few.
Ref. \cite{WaCh02a} investigated the effects of selective and random pinning schemes on scale-free networks.
Ref. \cite{ZhLuLu08}, by using adaptive controllers, provided a formula for estimating the number of pinned nodes and the strengths of couplings in a given complex network.
Ref. \cite{YuChLu09} found that a network can synchronize subject to any linear feedback pinning scheme by adaptively tuning the coupling strength.
Ref. \cite{SodiGaCh07,WaDaDoCaSu08} studied pinning controllability of complex networks in terms of the spectral properties of networks.
Refs. \cite{WaWuHuReWu16a,ChYuSu17,YuChWaYa09,ZhCaTaXiPe18,WaWeWuHuXu18,WeYuChYuCh18,WuXuPaShWa18} discussed pinning synchronization in certain networks, such as sensor networks, neural networks, multilayer networks, and multi-agent systems.
However, the above works mainly focused on proposing criteria for pinning synchronization, but did not investigate in depth how to choose the pinned nodes for optimal control.
In a recent survey \cite{waSu14} on pinning control, it shows that: there is a need for ``a deep understanding of the relationship between the effectiveness of pinning control and the structure properties of a network in the next decade'', and this understanding may facilitate better control of large-scale networked systems.

This paper revisits pinning control criteria for synchronization of complex dynamical networks, and proposes that
the effectiveness of a node-selection scheme is measured by the smallest eigenvalue ($\lambda_1$) of the submatrix obtained by deleting the rows and columns corresponding to the pinned nodes from the Laplacian matrix of the network. That submatrix is also referred to a grounded Laplacian matrix \cite{Mi93,BaHe06}. The larger the $\lambda_1$, the more effective the pinning control. Two key questions about the node-selection scheme emerging from the measure $\lambda_1$ are: to achieve an optimal performance, how many nodes does one need to control and which nodes should one apply controllers to? To answer these questions is by no means easy. Even when the number of pinned nodes is given, it is almost impossible to find the maximum of $\lambda_1$ by numerical calculations, especially in a large network. That is because it is an NP-hard problem with a high computational complexity to calculate the eigenvalues for all possible matrices obtained by node combinations in a network.

Grounded Laplacian matrix has attracted increasing attention in recent years. In \cite{PiSu16,PiSu14,PiShSu15}, upper and lower bounds on the smallest eigenvalue of grounded Laplacian matrix are derived, based on the sum of the weights of the edges connected to the grounded nodes and the eigenvector associated with the eigenvalue.
In this paper, utilizing tools from graph theory and matrix analysis \cite{Ba10,HoJo94,Ha95}, spectral properties and several other bounds for $\lambda_1$ are analyze at a relatively high accuracy based on the network topology information, including the spectrum of the network Laplacian matrix, the minimal degree of uncontrolled nodes, the number of edges between the controlled node set and the uncontrolled node set, etc.
Furthermore, the new theoretical results are applied to optimizing the pinning control scheme for network synchronization.
Hence, it provides guidance for choosing a proper pinning control scheme, which is demonstrated by illustrative examples of different types of complex networks.
In particular, node-selection schemes in some typical networks, such as stars, double-star graphs, fully-connected graphs, are investigated.
Finally, some established theorems and proposed algorithms are applied to scale-free, small-world, and real-world networks.

The rest of the paper is organized as follows. In Section II, pinning synchronization criteria are revisited for complex dynamical networks, and based on an algebraic graph-theoretic condition, the measure $\lambda_1$ is proposed. In Section III, spectral properties of $\lambda_1$ are investigated, which provide effective guidance for the selection of pinned nodes. Section IV shows applications of the theoretical results to scale-free, small-world, and real-world networks. Finally, conclusions are drawn in Section V.

\section{Pinning Control Network Model and Its Graph-Theoretic Synchronization Criterion}\label{se:model}

Consider a controlled network of $N$ identical dynamical nodes with a diffusive coupling, described by
\begin{equation}\label{eq:model}
\dot{x}_i=f(x_i)-c\sum_{j=1}^{N}l_{ij}\,P\,x_j+u_i(x_1,\dots,x_N)\,,
\end{equation}
where $i=1,2,\dots,N$. In (\ref{eq:model}), vector $x_i\in \R^{n}$ is the state of node $i$, function $f(\cdot)$ describes the self-dynamics of node $i$, positive constant $c$ denotes the coupling strength of the network, $u_i$ is the controller applied at node $i$ and is to be designed, and inner coupling matrix $P:\R^n\to \R^n$ is positive semi-definite. The network topology can be conveniently described by a graph $\bb G = (\mathcal V, \mathcal E)$ with the node set $\mathcal{V} = \{1, \dots, N\}$ and edge set $\mathcal E \subseteq \mathcal V \times \mathcal V$. Matrix  $L_N=[l_{ij}]_{N\times N}$ is the Laplacian matrix of graph $\bb G$, defined as follows: if there is an edge from node $j$ to node $i$ ($i\neq j$), then $l_{ij}=-1$; $l_{ii}=-\sum_{j\neq i}\,l_{ij}$. Matrix $L_N$ is generally asymmetrical when the network is directed, but is symmetrical when the network is undirected.

Suppose that the target state of the network is $s(t)$, which satisfies
\begin{equation}\label{eq:s}
{\dot s}(t) = f(s(t))\,, \qquad s(0)=s_0 \,.
\end{equation}
The goal of pinning control is to select a part of the nodes from network (\ref{eq:model}) to control, such that the states of all the nodes in the network can synchronize to the target state $s(t)$.

\begin{assumption}\label{assu:cond_f}
There exists a sufficiently large positive constant $\alpha$, such that
$$(y-z)^{\top}\,\left[\left(f(y)-f(z)\right)-\alpha\,P(y-z)\right] \leq -\mu \|y-z\|^2$$
holds for some positive constant $\mu$ and for all vectors $y,z \in \R^n$.
\end{assumption}

Here, $\alpha$ is determined by $f(\cdot)$ and the inner coupling matrix $P$.
Assumption \ref{assu:cond_f} is a standard assumption on the node dynamics \cite{BeBeHa04,ChLuWu11,LiCaWu14,LiCaWuLuTs15,Wu07}, which has been used to study \emph{global} synchronization of a network. The physical meaning of this assumption is that the self-dynamics $f(\cdot)$ is not ``too nonlinear''.

\begin{lemma}[\cite{HoJo94}, Theorem 1.1.6]\label{lem:eigenvalue_polynomial}
Let $p(\lambda)$ be a given polynomial of degree $k$. If ($\lambda$, $x$) is an eigenvalue-eigenvector pair of matrix $A$, then ($p(\lambda)$, $x$) is an eigenvalue-eigenvector pair of $p(A)$. Conversely, if $k\geq 1$ and if $\mu$ is an eigenvalue of $p(A)$, then there is some eigenvalue $\lambda$
of $A$ such that $\mu = p(\lambda)$.
\end{lemma}

For a real and symmetric matrix $A$, notation $A \succ 0$ means $A$ is positive definite.
\begin{lemma}[\textit{Schur Complement} \cite{BoElFeBa94}]\label{lem:schur}
For the symmetric matrix
$\begin{bmatrix}
      A & B\\
      B^{\top} & C
    \end{bmatrix}$,
where $A^{\top}=A$ and $C^{\top}=C$, the following three conditions are equivalent:\\
{\textit 1)} $\begin{bmatrix}
      A & B\\
      B^{\top} & C
    \end{bmatrix}\succ 0\,;$\\
{\textit 2)} $A\succ {0}$ and $C\,-\,B^{\top}\,A^{-1}\, B\succ {0}\,;$\\
{\textit 3)} $C\succ {0}$ and $A\,-\,B\,C^{-1}\, B^{\top}\succ {0}\,.$
\end{lemma}

\begin{lemma}\label{lem:from_shur}
Suppose $G=\begin{bmatrix}G_1& G_3 \\ G_3^{\top} & G_2  \end{bmatrix}$ and $D=\begin{bmatrix} D_1 & \bf{0} \\ \bf{0} & \bf{0} \end{bmatrix}$,
where $G, D \in \R^{N\times N}$, $G_1, D_1 \in \R^{l\times l}$ ($1\leq l <N$), $D_1 = \diag\{d,\dots,d\}$ with $d>0$, and $G_1$ and $G_2$ are symmetric. If $G_2\succ0$ and $d>\lambda_{\max}(G_3 G_2^{-1} G_3^{\top}-G_1)$, then $G+D\succ 0$.
\end{lemma}

\begin{IEEEproof}
For given matrices $G_1,\,G_2,\,G_3$ and $G_2\succ 0$, one can choose a sufficiently large constant $d$ such that $G_1+D_1-G_3 G_2^{-1} G_3^{\top}\succ 0$. It then follows from Lemma \ref{lem:schur} that
\begin{equation*}
G+D=\begin{bmatrix}
G_1+D_1&G_3\\
G_3^{\top}&G_2
\end{bmatrix}\succ 0\,.
\end{equation*}
Note that $G_1+D_1-G_3 G_2^{-1} G_3^{\top}\succ0$ is equivalent to $\lambda_{\min}(G_1+dI_l-G_3 G_2^{-1} G_3^{\top})>0$. In addition, by Lemma \ref{lem:eigenvalue_polynomial}, one has $\lambda_{\min}(G_1+dI_l-G_3 G_2^{-1} G_3^{\top})=\lambda_{\min}(G_1-G_3 G_2^{-1} G_3^{\top})+d$. Thus, $G_1+D_1-G_3 G_2^{-1} G_3^{\top}\succ0$ is equivalent to $d>\lambda_{\max}(G_3 G_2^{-1} G_3^{\top}-G_1)$.
\end{IEEEproof}

Arrange all eigenvalues of a real and symmetric matrix $M\in \R^{N\times N}$ in a non-increasing order as follows:
\begin{equation}\label{eq:eigenvalue_rank}
\lambda_N(M)\geq\lambda_{N-1}(M)\geq\dots\geq\lambda_1(M)\,.
\end{equation}
Here, $\lambda_1(M)$ is the smallest eigenvalue of $M$ and $\lambda_N(M)$ is the biggest one.

Let $\mathcal S = \{s_1,s_2,\dots,s_l\}$, where $\mathcal S \subset \mathcal V$ and
$l < N$, be the set of controlled nodes in the network $\bb G$. Without loss of generality, re-order the nodes in the network to let the first $l$ ($1\leq l < N$) nodes be controlled. Let $\hat{L}_N$ be the symmetrized matrix of $L_N$ and $\hat{L}_N=(L_N+L^{\top}_N)/2$. The symmetrized Laplacian matrix $\hat{L}_N$ can be written in the block form of $\begin{bmatrix}
\hat{L}_l&\hat{L}_{12}\\
\hat{L}_{21}&\hat{L}_{N-l}
\end{bmatrix}$.
In this paper, $I_{N}$ denotes the identity matrix of order $N$.

Let $e_i(t)\,\triangleq\,x_i(t)-s(t)$. From the controlled network (\ref{eq:model}) and the target dynamics (\ref{eq:s}), the error system is obtained as
\begin{equation}\label{eq:e_i}
    \dot{e}_i=f(x_i)-f(s)-c\sum\limits_{j=1}^N\,l_{ij}\,P\,e_j+u_i(x_1,\dots,x_N)\,,
\end{equation}
where $i=1,\dots,N$.

The objective now is to design controllers $u_i$ such that (\ref{eq:model}) and (\ref{eq:s}) together satisfy
\begin{equation}
\lim_{t \rightarrow \infty}\|x_i(t)-s(t)\| = 0\,,  \quad \text{for}~i=1,2,\dots,N\,.
\end{equation}
Design $u_i$ by using two types of controllers: adaptive controllers and linear feedback ones.
First, design adaptive pinning controllers of the form
\begin{equation}\label{eq:adaptive_con}
\begin{cases}
u_{i} = -{d_iPe_i}\,,\quad \dot{d_i}=h_ie^{\top}_iPe_i\,,& 1\leq i\leq l\,,\\
u_i=0\,,                                                 & l+1\leq i\leq N\,,
\end{cases}
\end{equation}
where $h_i$, $1\leq i\leq l$, are positive constants.

\begin{theorem}[\cite{SoCa10,YuChWaYa09,WaDaDoCaSu08}]\label{the:model}
Suppose that Assumption \ref{assu:cond_f} holds. The states of all the nodes in network (\ref{eq:model}) synchronize to the target state $s(t)$ described by (\ref{eq:s}) using the adaptive pinning controller (\ref{eq:adaptive_con}),
provided that
\begin{equation}\label{eq:con}
\lambda_{1}({\hat L}_{N-l})>\alpha/c \,.
\end{equation}
Here, ${\hat{L}}_{N-l}$ is the principal submatrix of $\hat{L}_N$ obtained by deleting the first $l$ rows and columns from $\hat{L}_N$.
\end{theorem}

\begin{IEEEproof}
Consider the Lyapunov function
$$V(t)=\frac{1}{2}\sum^N_{i=1} e^{\top}_i e_i + \frac{1}{2}\sum^l_{i=1} \frac{(d_i - d^*)^2}{h_i}\,,$$
where $d^*$ is a sufficiently large positive constant to be determined later in the proof. Calculating the derivative of $V(t)$ using (\ref{eq:e_i}) and (\ref{eq:adaptive_con}) yields
\begin{equation*}
\begin{split}
\frac{dV}{dt}=&\sum^N_{i=1}e^{\top}_i \dot{e}_i + \sum^l_{i=1}\frac{\dot{d_i}(d_i-d^*)}{h_i} \\
             =&\sum^N_{i=1}e^{\top}_i [f(x_i)-f(s)-\alpha P(x_i-s)]\\
              &+\alpha \sum^N_{i=1}e^{\top}_iPe_i-c\sum^N_{i=1}e^{\top}_i\sum^N_{j=1}l_{ij}Pe_j\\
              &-\sum^l_{i=1}d_i e^{\top}_iPe_i+\sum^l_{i=1}(d_i-d^*)e^{\top}_iPe_i\,.
\end{split}
\end{equation*}
By Assumption \ref{assu:cond_f}, one has
\begin{equation*}
\begin{split}
\frac{dV}{dt}\leq&-\mu\sum^N_{i=1}\|e_i\|^2 + \alpha \sum^N_{i=1}e^{\top}_iPe_i - d^* \sum^l_{i=1}e^{\top}_iPe_i\\
              &- c\sum^N_{i=1}e^{\top}_i\sum^N_{j=1}l_{ij}Pe_j\,.
\end{split}
\end{equation*}
Let $e=[e^{\top}_1,e^{\top}_2,\dots,e^{\top}_N]^{\top}\in \R^{nN}$ and $\otimes$ denote the Kronecker product. The above inequality can be rewritten as
\begin{equation}\label{eq:19}
\begin{split}
\frac{dV}{dt}\leq &-\mu e^{\top}e + \alpha e^{\top}(I_N\otimes P)e \\
              &- d^*e^{\top}
\Bigg(\begin{bmatrix}
     I_l&\\
     &{\bf{0}}_{N-l}
\end{bmatrix}\otimes P\Bigg)e - ce^{\top}({L}_N\otimes P)e \\
          %\leq&\,e^{\top}[(\alpha I_N -d
%\begin{bmatrix}
%     I_l&\\
%     &{\bf{0}}_{N-l}
%\end{bmatrix}-c\hat{L}_N)\otimes P]e\\
             =&-\mu e^{\top}e-e^{\top}[(-\alpha I_N +d^*
\begin{bmatrix}
     I_l&\\
     &{\bf{0}}_{N-l}
\end{bmatrix}+c\hat{L}_N)\otimes P]e \\
=&-\mu e^{\top}e-e^{\top}[(-\alpha I_N +D^*
+c\hat{L}_N)\otimes P]e\,, \\
\end{split}
\end{equation}
where $D^*\triangleq d^*\begin{bmatrix}
     I_l&\\
     &{\bf{0}}_{N-l}
\end{bmatrix}\,.$
Note that $$-\alpha I_N +D^*+c\hat{L}_N = \begin{bmatrix}
c\hat{L}_l-\alpha I_l &c\hat{L}_{12}\\
c\hat{L}_{21}&c\hat{L}_{N-l}-\alpha I_{N-l}
\end{bmatrix}+D^*\,.$$
If the symmetric matrix $c\hat{L}_{N-l}-\alpha I_{N-l}\succ 0$, one can choose a sufficiently large constant $d^*$ such that $-\alpha I_N +D^* +c\hat{L}_N\succ 0$, according to Lemma \ref{lem:from_shur}.
Noticing that $P$ is positive semi-definite in network (\ref{eq:model}), one has  $(-\alpha I_N +D^* +c\hat{L}_N)\otimes P \succeq 0$. Therefore, (\ref{eq:19}) reduces to $\frac{dV}{dt}\leq -\mu e^{\top}e$. Thus, $\frac{dV}{dt}\leq 0$ and the  equality holds if and only if $e=\bf{0}$. According to the Lyapunov stability theorem \cite{Kh02}, $e_i$ in system (\ref{eq:e_i}) converges to $\bf{0}$ as $t \rightarrow \infty$, for all $i=1,2,\dots,N$.\\
Finally, it follows from Lemma \ref{lem:eigenvalue_polynomial} that the condition $c\hat{L}_{N-l}-\alpha I_{N-l}\succ 0$ is equivalent to $\lambda_1(\hat{L}_{N-l})>\alpha/c$. The proof is completed.
\end{IEEEproof}

Next, pinning control of network (\ref{eq:model}) is designed by using linear feedback controllers. Suppose that the first $l$ ($1\leq l < N$) nodes are controlled under the linear state-feedback controllers \cite{Ch14} of the form
\begin{equation}\label{eq:linear_con}
\begin{cases}
u_i=-cdPe_i\,, & 1\leq i\leq l\,,\\
u_i=0\,, & l+1\leq i\leq N\,,
\end{cases}
\end{equation}
where $d$ is a positive constant feedback gain for all the $l$ controllers and it is to be determined.

\begin{theorem}[\cite{SoCa10,YuChWaYa09}]\label{the:linear model}
Suppose that Assumption \ref{assu:cond_f} holds. The states of network (\ref{eq:model}) synchronize to the target state $s(t)$ by using the linear feedback controller (\ref{eq:linear_con}),
provided that
$\lambda_{1}({\hat L}_{N-l})>\alpha/c$
and the feedback gain satisfies
\begin{equation}\label{eq:lin_con_2}
d>\frac{1}{c}(\lambda_{\max}(c^2\hat{L}_{12} (c\hat{L}_{N-l}-\alpha I_{N-l})^{-1} \hat{L}_{21}-c\hat{L}_l )+\alpha)\,.
\end{equation}
\end{theorem}

\begin{IEEEproof}
Let the gain matrix $D \triangleq d\begin{bmatrix}
I_l&\\
&{\bf{0}}_{N-l}
\end{bmatrix}_{N\times N}$,
and consider a Lyapunov function $V(t)=\frac{1}{2}\sum^{N}_{i=1}e^{\top}_{i}e_{i}\,.$ Using (\ref{eq:e_i}) and (\ref{eq:linear_con}), the derivative of $V(t)$ gives
\begin{equation*}
\begin{split}
\frac{dV}{dt}
            %=&\sum^N_{i=1}e^{\top}_i [f(x_i)-f(s)-c\sum^N_{j=1}l_{ij}Pe_j]+\sum^l_{i=1}e^{\top}_i u_i\\
             =&\sum^N_{i=1}e^{\top}_i [f(x_i)-f(s)-\alpha P(x_i-s)]+\alpha\sum^N_{i=1}e^{\top}_iPe_i\\
              &-c\sum^N_{i=1}e^{\top}_i\sum^N_{j=1}l_{ij}Pe_j-c\sum^l_{i=1}d e^{\top}_iPe_i\\
\leq  &   -\mu e^{\top}e + \alpha e^{\top}(I_N\otimes P)e \\
              &- ce^{\top}(D\otimes P)e - ce^{\top}(L_N\otimes P)e \\
          =& -\mu e^{\top}e
             -e^{\top}[(-\alpha I_N +cD+c\hat{L}_N)\otimes P]e\,.
\end{split}
\end{equation*}
Note that
\begin{equation*}
\begin{split}
&-\alpha I_N+cD +c\hat{L}_N\\
&=\begin{bmatrix}
c\hat{L}_l-\alpha I_l &c\hat{L}_{12}\\
c\hat{L}_{21}&c\hat{L}_{N-l}-\alpha I_{N-l}
\end{bmatrix}+cd\begin{bmatrix}
I_l&\\
&{\bf{0}}_{N-l}
\end{bmatrix}\,.
\end{split}
\end{equation*}
By Lemma \ref{lem:from_shur}, if
%\begin{equation}\label{eq:7}
$c\hat{L}_{N-l}-\alpha I_{N-l}\succ 0$
%\end{equation}
and
\begin{equation}\label{eq:8}
cd>\lambda_{\max}(c^2\hat{L}_{12} (c\hat{L}_{N-l}-\alpha I_{N-l})^{-1} \hat{L}_{21}-(c\hat{L}_l-\alpha I_l))\,,
\end{equation}
then one has $-\alpha I_N+cD +c\hat{L}_N\succ 0$. It follows that $(-\alpha I_N +D +c\hat{L}_N)\otimes P \succeq 0$ and $\frac{dV}{dt}\leq -\mu e^{\top}e$. Therefore, $\frac{dV}{dt}\leq 0$ and the equality holds if and only if $e=\bf{0}$. Then, $e_i$ in system (\ref{eq:e_i}) converges to $\bf{0}$ as $t \rightarrow \infty$ for all $i=1,2,\dots,N$.
Finally, it is easy to verify that (\ref{eq:8}) can be reduced to (\ref{eq:lin_con_2}).
This completes the proof.
\end{IEEEproof}

%\begin{remark}
A wide spectrum of research on pinning control of complex dynamical networks has been performed in the last few years, such as in \cite{LiWaCh04,ZhLuLu08,Wu08a,YuChLu09,SoCa10,Ch14}.
Some closely related results to Theorems \ref{the:model} \& \ref{the:linear model} are compared and discussed as follows.
%\end{remark}

\begin{remark}
Similar graph-theoretic criteria were proposed in \cite{SoCa10,YuChWaYa09,WaDaDoCaSu08}. However, the assumption on self-dynamics of individual node therein is different from Assumption \ref{assu:cond_f} here.
%Thus, the corresponding proofs are different to some extent.
The requirement on individual dynamics in Assumption \ref{assu:cond_f} is weaker than that in \cite{SoCa10,YuChWaYa09,WaDaDoCaSu08}, since different inner-coupling patterns can be assigned to $P$ through design.
%In order to maintain the integrity of the paper, short proofs for Theorems \ref{the:model} and \ref{the:linear model} are kept in this section.
The proofs for Theorems \ref{the:model} and \ref{the:linear model} are similar to those in \cite{SoCa10,YuChWaYa09,WaDaDoCaSu08}, but different to some extent due to the weaker Assumption \ref{assu:cond_f}.
%Interested readers can refer to our online note \cite{LiXuChLu18} for the detailed proofs.
\end{remark}

\begin{remark}
Reference \cite{SodiGaCh07} defined the pinning controllability of a network in terms of the spectrum of an extended network topology.
The measure in \cite{SodiGaCh07} used the eigenratio or eigenvalues of the $(N+1)$-dimensional Laplacian matrix of the extended network topology, and the Laplacian matrix is parameterized by the pinning control gains. In comparison, the measure introduced here uses the spectrum of the $(N-l)$-dimensional grounded matrix that is independent of the pinning control gains.
%The present work focuses on how to optimize the selection of pinned nodes resulting in a larger $\lambda_1$ of the grounded matrix under a given number of pinned nodes.
\end{remark}

From Theorems \ref{the:model} and \ref{the:linear model}, $\lambda_{1}({\hat L}_{N-l})>\alpha/c$ is the common algebraic graph-theoretic condition for pinning synchronization under the adaptive controllers or the linear state-feedback controllers.
%If the simple linear state-feedback controllers are used, one more condition is that the feedback gain $d$ needs to be large enough to satisfy (\ref{eq:lin_con_2}).
In the rest of this paper, only undirected networks are discussed. In this case, the pinning synchronization criterion is simplified to be $\lambda_{1}({L}_{N-l})>\alpha/c$, because $\hat{L}_N=L_N$.

\section{Effective Pinning Control Using Network Topology Information}\label{se:main}

From the above section, one can see that $c\lambda_1({L}_{N-l})>\alpha$ determines the pinning-control synchronizability. Under Assumption \ref{assu:cond_f}, $\alpha$ is  determined by the given node dynamics $f$ and the inner coupling pattern $P$. So, $\lambda_1(L_{N-l})$ can be used to measure the effectiveness of various pinning schemes. To be specific, for a given network, if the number of controlled nodes is limited or fixed beforehand, it is important to come up with a strategy to select a suitable set of pinned nodes from the network to maximize $\lambda_1({ L}_{N-l})$. In this section, tools from graph theory and matrix analysis \cite{Ba10,HoJo94,Ha95} will be used
to analyze the spectral properties of $\lambda_1({ L}_{N-l})$ based on the network topology information, so as to provide insights for the selection of pinned nodes.

\subsection{Notations and Useful Lemmas}

Let $\mathcal S = \{s_1,s_2,\dots,s_l\}$ be the set of controlled nodes, and $|\mathcal S|$ denote the number of nodes in set $\mathcal S$, so $|\mathcal S| = l$. Following the notations used in \cite{Ba10}, $L(\mathcal S|\mathcal S)$ is the $(N-l)\times(N-l)$ principal submatrix of $L_N$ obtained by deleting the rows and columns corresponding to $\mathcal S$. In the rest of this paper, $L_{N-l}$ can be understood as $L(\mathcal S|\mathcal S)$. Specially, when $\mathcal S = \{s_1\}$ is a singleton, $L(\mathcal S|\mathcal S)$ is denoted as $L(s_1|s_1)$.

Note that the sufficient condition (\ref{eq:con}) for the synchronization in an undirected network can be written as
\begin{equation}
\lambda_1(L(\mathcal S|\mathcal S)) > \alpha/c, \quad \text{where}~\mathcal S = \{s_1,s_2,\dots,s_l\}\,.
\end{equation}

The following lemmas from matrix analysis are useful for analyzing properties of $\lambda_{1}(L(\mathcal S|\mathcal S))$.

%\begin{lemma}(Spielman, 2004; Anderson {\em \&} Morley, 1985)\label{lem:laplacian_min_0}
%If graph $\bb G$ has $N$ nodes and is undirected, and $\lambda$ is an eigenvalue of $L(\bb{G})$, then $0\leq \lambda \leq N$.  The multiplicity of $0$ equals the number of the connected components of $\bb G$.
%\end{lemma}

\begin{lemma}[\cite{HoJo94}, Theorem 1.1.6]\label{lem:eigenvalue_polynomial}
Let $p(\lambda)$ be a given polynomial of degree $k$. If ($\lambda$, $x$) is an eigenvalue-eigenvector pair of matrix $A$, then ($p(\lambda)$, $x$) is an eigenvalue-eigenvector pair of $p(A)$. Conversely, if $k\geq 1$ and if $\mu$ is an eigenvalue of $p(A)$, then there is some eigenvalue $\lambda$
of $A$ such that $\mu = p(\lambda)$.
\end{lemma}

\begin{lemma}[\cite{Ha95}, Cauchy Interlace Theorem]\label{lem:eig_com_0}
Let $U$ be a real $n \times m$ matrix such that $U^{\top}U = I_m$ and let $A$ be a symmetric $n \times n$ matrix with eigenvalues $\lambda_n\geq\dots\geq\lambda_1$. Define $B = U^{\top}AU$ and let $B$ have eigenvalues $\mu_{m}\geq\dots\geq\mu_1$. Then, the eigenvalues of $B$ interlace those of $A$, i.e.,
\begin{equation*}
    \lambda_{n-m+i}\geq\mu_{i}\geq\lambda_{i}\,.
\end{equation*}
In addition, if $C$ is a principal submatrix of a symmetric matrix $A$, then the eigenvalues of $C$ interlace the eigenvalues of $A$.
\end{lemma}

\begin{lemma}[\cite{Ba10}, Chapter 1]\label{lem:eig_com}
Let $A$ be a real symmetric $N\times N$ matrix and $B$ be a principal submatrix of $A$ of order $N-1$. If $\lambda_N\geq\dots\geq\lambda_1$ and $\mu_{N-1}\geq\dots\geq\mu_1$ are the eigenvalues of $A$ and $B$, respectively, then
\begin{equation*}
    \lambda_N\geq\mu_{N-1}\geq\lambda_{N-1}\geq\dots\geq\lambda_2\geq\mu_1\geq\lambda_1\,.
\end{equation*}
\end{lemma}

Lemma \ref{lem:eig_com} is a special case of Lemma \ref{lem:eig_com_0}.

\begin{lemma}[\cite{Ba10}, Chapter 1]\label{lem:eigenvalue_definition}
Let $A$ be a real symmetric $N\times N$ matrix with eigenvalues $\lambda_N(A)\geq\dots\geq\lambda_1(A)$. Let $\|x\| = (\sum_{i=1}^{N}x^2_i)^{1/2}$. Then,
\begin{equation}\label{eq:10}
    \lambda_N(A)=\max_{\|x\|=1}[x^{\top}Ax],\quad\lambda_1(A)=\min_{\|x\|=1}[x^{\top}Ax]\,.
\end{equation}
\end{lemma}

\begin{lemma}[\cite{Ba10}, Chapter 1]\label{lem:eig_diag}
Let $A=[a_{ij}]$ be a real symmetric $N\times N$ matrix. Then,
\begin{equation*}
    \lambda_N(A) \geq \max_{1\leq i\leq N}[a_{ii}] \geq \min_{1\leq i\leq N}[a_{ii}] \geq \lambda_1(A).
\end{equation*}
\end{lemma}

Lemma \ref{lem:eig_diag} shows that the smallest eigenvalue of a real symmetric matrix is no greater than the smallest diagonal element of the matrix.

\begin{lemma}[\cite{HoJo94}, Theorem 4.3.1]\label{lem:weyl}
Let real matrices $A,B$ be symmetric and let the respective eigenvalues of $A,B$, and $A+B$ be $\{\lambda_i(A)\}_{i=1}^N$, $\{\lambda_i(B)\}_{i=1}^N$ and $\{\lambda_i(A+B)\}_{i=1}^N$, each ordered as in (\ref{eq:eigenvalue_rank}). Then,
\begin{equation}\label{eq:13}
\lambda_i(A+B)\leq \lambda_{i+j}(A)+\lambda_{N-j}(B),\quad j=0,1,\dots,N-i
\end{equation}
for each $i=1,\dots,N$.
Also,
\begin{equation}\label{eq:14}
\lambda_{i-j+1}(A)+\lambda_{j}(B)\leq \lambda_i(A+B),\quad j=1,\dots,i
\end{equation}
for each $i=1,\dots,N$.
\end{lemma}

\begin{definition}[\cite{Di00}]
For graph $\bb{G} = (\mathcal{V}, \mathcal{E})$, if $\mathcal{V}^{'} \subseteq \mathcal V$ and $\mathcal{E}^{'} \subseteq \mathcal E$, then $\bb{G}^{'} = (\mathcal{V}^{'}, \mathcal{E}^{'})$ is a \emph{subgraph} of $\bb G$. Furthermore, if $\bb{G}^{'}$ contains all the edges $(i, j) \in \mathcal E$ with $i, j \in \mathcal{V}^{'}$, then $\bb{G}^{'}$ is an \emph{induced subgraph} of $\bb G$.
\end{definition}

In graph $\bb G$, let $k_i$ denote the degree of node $i$.
Let $p_1,\dots,p_{N-l}$ denote the nodes in set $\mathcal{V}\setminus\mathcal S$. Let $w_{p_j}$ be the number of nodes in $\mathcal S$ that are connected to node $p_j$, so $w_{p_j} \geq 0$.  Let $\bb H$ be the induced subgraph of $\bb G$ with node set $\mathcal V\setminus\mathcal S$. Then, matrix $L(\mathcal S|\mathcal S)$ can be written as
\begin{equation}\label{eq:trans}
L(\mathcal S|\mathcal S) = L(\bb H)+\Lambda\,,
\end{equation}
where $L(\bb H)$ is the Laplacian matrix of graph $\bb H$ and $\Lambda=\diag\{w_{p_1},\dots,w_{p_{N-l}}\}$.

%\begin{definition} (Horn {\em \&} Johnson, 1994, Definition 6.2.25)
%Let $\mathbb C$ be the set of complex numbers and matrix $A\in \mathbb{C}^{n\times n}$. Matrix $A$ is \emph{irreducibly diagonally dominant} if\\
%(a) $A$ is irreducible;\\
%(b) $A$ is diagonally dominant, that is, $|a_{ii}|\geq R^{'}_{i}(A)$ for all $i=1,\dots,n$, where $R^{'}_i(A)=\sum_{j=1,j\neq i}^{n}|a_{ij}|$;\\
%(c) There is an $i\in \{1,\dots,n\}$ such that $|a_{ii}|> R^{'}_{i}(A)$.
%\end{definition}
%
%\begin{lemma}\label{lem:irreducible}
%(Horn {\em \&} Johnson, 1994, Corollary 6.2.27)
%Let $A\in \mathbb{C}^{n\times n}$ be irreducibly diagonally dominant.
%%(a) $A$ is nonsingular;\\
%%(b) if every main diagonal entry of $A$ is real and positive, then every eigenvalue of $A$ has positive real part;\\
%If $A$ is Hermitian and every main diagonal entry is positive, then every eigenvalue of $A$ is positive, that is, $A$ is positive definite.
%\end{lemma}
%
%Lemma \ref{lem:irreducible} is useful in the case that matrix $A$ is real and %symmetric.

\subsection{Spectral Properties and Their Applications to Optimizing Pinning Control}

\begin{theorem}[\cite{HoJo94}]\label{the:eigenpositive}
Suppose that network $\bb G$ is undirected and connected. If the set of controlled nodes $\mathcal S$ is not empty,
then the matrix $L(\mathcal S|\mathcal S)$ is positive definite, i.e., $\lambda_{1}(L(\mathcal S|\mathcal S))>0$ for $1\leq l\leq N-1$.
\end{theorem}

Theorem \ref{the:eigenpositive} is easy to verify. Note that every irreducible component of $L(\mathcal S|\mathcal S)$ is diagonally dominant with at least one strictly diagonally dominant row. Thus, $L(\mathcal S|\mathcal S)$ is positive definite according to Corollary 6.2.27 in \cite{HoJo94}.

\begin{remark}
Theorem \ref{the:eigenpositive} shows that no matter which and how many nodes are controlled, the smallest eigenvalue satisfies $\lambda_{1}(L(\mathcal S|\mathcal S))>0$ ($1\leq l\leq N-1$) for the undirected and connected network. It means that even only one node is pinned, a sufficiently large $c$ can make $c\lambda_{1}(L(\mathcal S|\mathcal S))>\alpha$, so that the network synchronizes.
This is in accordance with the conclusion in \cite{ChLilu07}.
\end{remark}

\begin{theorem}\label{the:mul_laplacian_spectrum}
Let $\mathcal S$ be the set of $l$ controlled nodes, which are arbitrarily chosen from $\mathcal V$.
Then,
\begin{equation}\label{eq:38}
    \lambda_{1}(L(\mathcal S|\mathcal S))\leq\lambda_{l+1}(L_N)\,.
\end{equation}
\end{theorem}
\begin{IEEEproof}
$L(\mathcal S|\mathcal S)$ is an $(N-l)\times(N-l)$ principal submatrix of $L_N$. By Lemma \ref{lem:eig_com_0}, one has
\begin{equation}\label{eq:40}
\lambda_{l+i}(L_N)\geq\lambda_i(L(\mathcal S|\mathcal S))\geq\lambda_{i}(L_N) \,.
\end{equation}
Let $i=1$ in (\ref{eq:40}). Then, $\lambda_{1}(L(\mathcal S|\mathcal S))\leq\lambda_{l+1}(L_N)$.
\end{IEEEproof}

%\begin{corollary}\label{co:mul_laplacian_spectrum}
%A necessary condition for $\lambda_1(L(\mathcal S|\mathcal S)) > \alpha/c$, where $|\mathcal S| = l$, is that
%\begin{equation}\label{eq:37}
%\lambda_{l+1}(L_{N}) > \alpha/c\,.
%\end{equation}
%\end{corollary}

From (\ref{eq:38}), one can predict the minimal number of nodes that must be pinned through analyzing the Laplacian spectrum of the network, to guarantee $\lambda_1(L(\mathcal S|\mathcal S)) > \alpha/c$.

Consider the situation where only one node,  denoted by $i$, is controlled. That is, $\mathcal S = \{i\}$. From Theorem \ref{the:mul_laplacian_spectrum}, one has $\lambda_1(L(i|i))\leq \lambda_2(L_N)$. Thus, the following result holds.

\begin{corollary}\label{co:one_laplacian_spectral}
A necessary condition for $\lambda_{1}({L}(i|i))>\alpha/c$ is that
\begin{equation}\label{eq:21}
   \lambda_{2}(L_N)>\alpha/c\,.
\end{equation}
\end{corollary}

\begin{remark}
Here, $c\cdot\lambda_2(L_N)>\alpha$ is the synchronization condition for a complex network of coupled oscillators without pinning control, see \cite{Wu03}. And $\lambda_2(L_N)$ measures the synchronizability of the network, see \cite{PeCa98,WaCh02b}. In addition, $\lambda_2(L_N) > 0$ is a necessary condition for $\lambda_1(L(i|i)) > 0$, which implies that the network has to be connected in order to synchronize the network by pinning only one node.
\end{remark}

Next, two examples are given to show the effectiveness of the estimation of $\lambda_1(L(\mathcal S|\mathcal S))$ in Theorem \ref{the:mul_laplacian_spectrum}.

\begin{exmp}\label{ex:two_stars}
Consider the double-star graph shown in Fig. \ref{fig:3}. The number of nodes is $N=13$. Its Laplacian matrix is
\begin{equation*}
L_{N}=\left[\begin{smallmatrix}
2& {-1} &0&0&0&0&0& {-1} &0&0&0&0&0\\
-1&6&-1&-1&-1&-1&-1&0&0&0&0&0&0\\
0&-1&1&0&0&0&0&0&0&0&0&0&0\\
0&-1&0&1&0&0&0&0&0&0&0&0&0\\
0&-1&0&0&1&0&0&0&0&0&0&0&0\\
0&-1&0&0&0&1&0&0&0&0&0&0&0\\
0&-1&0&0&0&0&1&0&0&0&0&0&0\\
-1&0&0&0&0&0&0&6&-1&-1&-1&-1&-1\\
0&0&0&0&0&0&0&-1&1&0&0&0&0\\
0&0&0&0&0&0&0&-1&0&1&0&0&0\\
0&0&0&0&0&0&0&-1&0&0&1&0&0\\
0&0&0&0&0&0&0&-1&0&0&0&1&0\\
0&0&0&0&0&0&0&-1&0&0&0&0&1
\end{smallmatrix}\right]\,.
\end{equation*}

\begin{figure}
\centering
\includegraphics[width=2.5in]{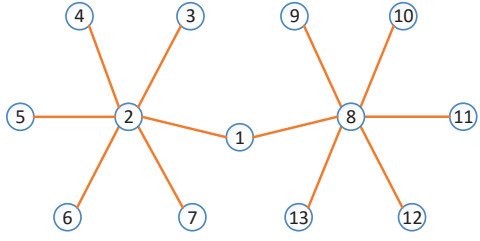}
\caption{A double-star graph.}
\label{fig:3}
\end{figure}

Arrange all the eigenvalues of $L_N$ in a non-decreasing order:
\begin{equation}\label{eq:Lspectra}
\begin{split}
\{\lambda_i(L_N),\,i=1,\dots,N~~|&~~0~~0.1459~~1~~1~~1~~1~~1~~1~~1~~1\\
                                 &~~1.8074~~6.8541~~7.1926\}\,.
\end{split}
\end{equation}
Then, calculate $\lambda_1(L(\mathcal S|\mathcal S))$ for $l=1,\dots,12$. Since there are $C_N^l$ options for choosing $l$ nodes from $N$ nodes, there are $C_N^l$ possible values for $\lambda_{1}({L} (\mathcal S|\mathcal S))$. Let $\max_{|{\mathcal S}|=l}[\lambda_{1}({L} (\mathcal S|\mathcal S))]$ be the maximum of those values. Then,
\begin{equation*}
\begin{split}
&\max_{|{\mathcal S}|=1}[\lambda_{1}(L(\mathcal S|\mathcal S))]=0.1459\,;\\
&\max_{|{\mathcal S}|=2}[\lambda_{1}(L(\mathcal S|\mathcal S))]=1,\dots,\max_{|{\mathcal S}|=9}[\lambda_{1}(L(\mathcal S|\mathcal S))]=1\,;\\
&\max_{|{\mathcal S}|=10}[\lambda_{1}(L(\mathcal S|\mathcal S))]=1.5505\,;\\
&\max_{|{\mathcal S}|=11}[\lambda_{1}(L(\mathcal S|\mathcal S))]=6\,;\\
&\max_{|{\mathcal S}|=12}[\lambda_{1}(L(\mathcal S|\mathcal S))]=6\,.
\end{split}
\end{equation*}

Figure \ref{fig:4} shows the trends of $\max_{|{\mathcal S}|=l}[\lambda_{1}(L(\mathcal S|\mathcal S))]$ and $\lambda_{l+1}(L_N)$ as $l$ increases from 1 to 12. From Fig. \ref{fig:4}, one can see that $\lambda_{l+1}(L_N)$ keeps being greater than or equal to $\lambda_1(L(\mathcal S|\mathcal S))$ for increasing $l$. It verifies Theorem \ref{the:mul_laplacian_spectrum}. Furthermore, $\lambda_{l+1}(L_N)=\max_{|\mathcal S|=l}[\lambda_1(L(\mathcal S|\mathcal S))]$ for $l=1,\dots,9$, which means that $\lambda_{l+1}(L_N)$ can well estimate $\max_{|\mathcal S|=l}[\lambda_1(L(\mathcal S|\mathcal S))]$ when a small part of nodes ($\frac{l}{N} \leq 50\%$) are controlled.

\begin{remark}
$\max_{|\mathcal S|=l}[\lambda_1(L(\mathcal S|\mathcal S))]$ remains to be 1 as $l$ increases from 2 to 9, i.e., the percentage of pinned nodes increases from $15.4\%$ to $69.2\%$. It means that $\max_{|\mathcal S|=l}[\lambda_1(L(\mathcal S|\mathcal S))]$ is not improved as the number of pinned nodes $l$ increases. Figures \ref{fig:5}(a)-(d) show four pinning schemes for the set $\mathcal S$ that yields $\max_{|{\mathcal S}|=l}[\lambda_{1}(L(\mathcal S|\mathcal S))]$, when $l=2,3,4,9$, respectively. Since all the $\lambda_1(L(\mathcal S|\mathcal S))$ in the four cases are 1, it is cost-effective to choose the pinning scheme with the least number of controlled nodes, i.e., the scheme shown in Fig. \ref{fig:5}(a).
\end{remark}

\begin{figure}
\centering
\includegraphics[width=3in]{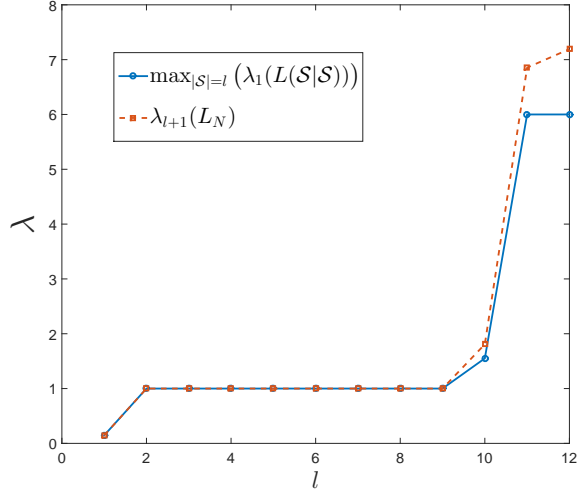}
\caption{The trends of $\max_{|{\mathcal S}|=l}[\lambda_{1}(L(\mathcal S|\mathcal S))]$ and $\lambda_{l+1}(L_N)$ for the graph in Fig. \ref{fig:3}, as the number of pinned nodes $l$ increases from 1 to 12.}
\label{fig:4}
\end{figure}

\begin{figure}
  \centering
  \subfigure[]{
    %\label{fig:subfig:a} %% label for first subfigure
    \includegraphics[width=1.5in]{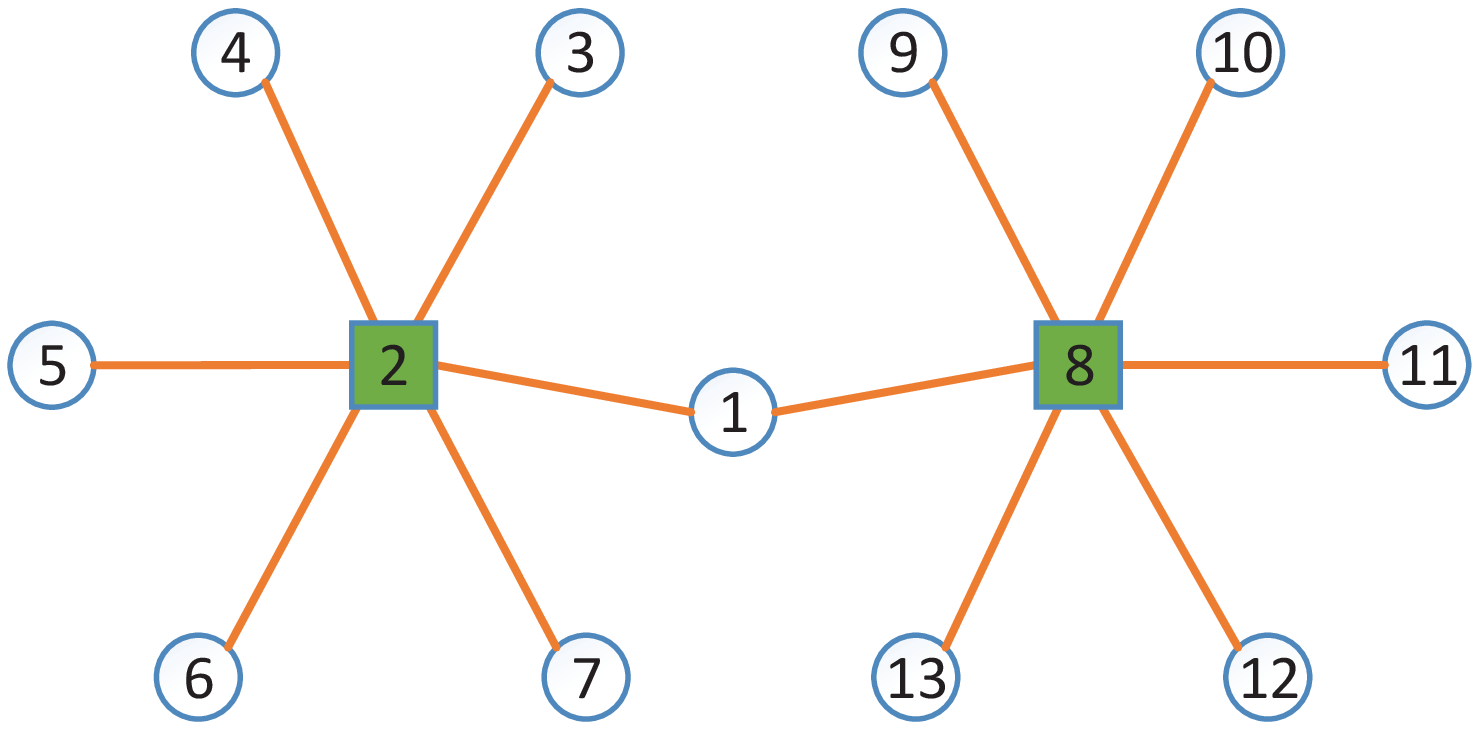}}
  \hspace{0.15in}
  \subfigure[]{
    %\label{fig:subfig:b} %% label for second subfigure
    \includegraphics[width=1.5in]{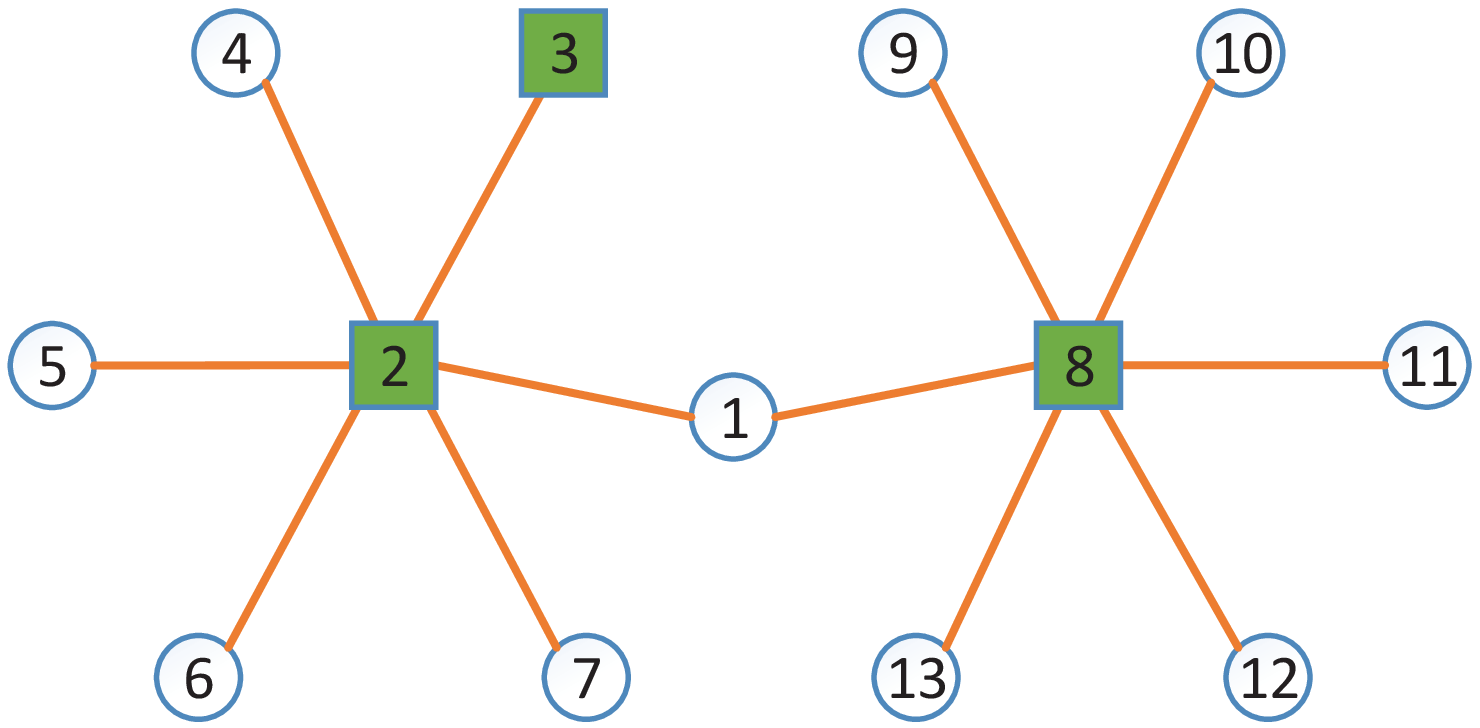}}
    \subfigure[]{
    %\label{fig:subfig:a} %% label for first subfigure
    \includegraphics[width=1.5in]{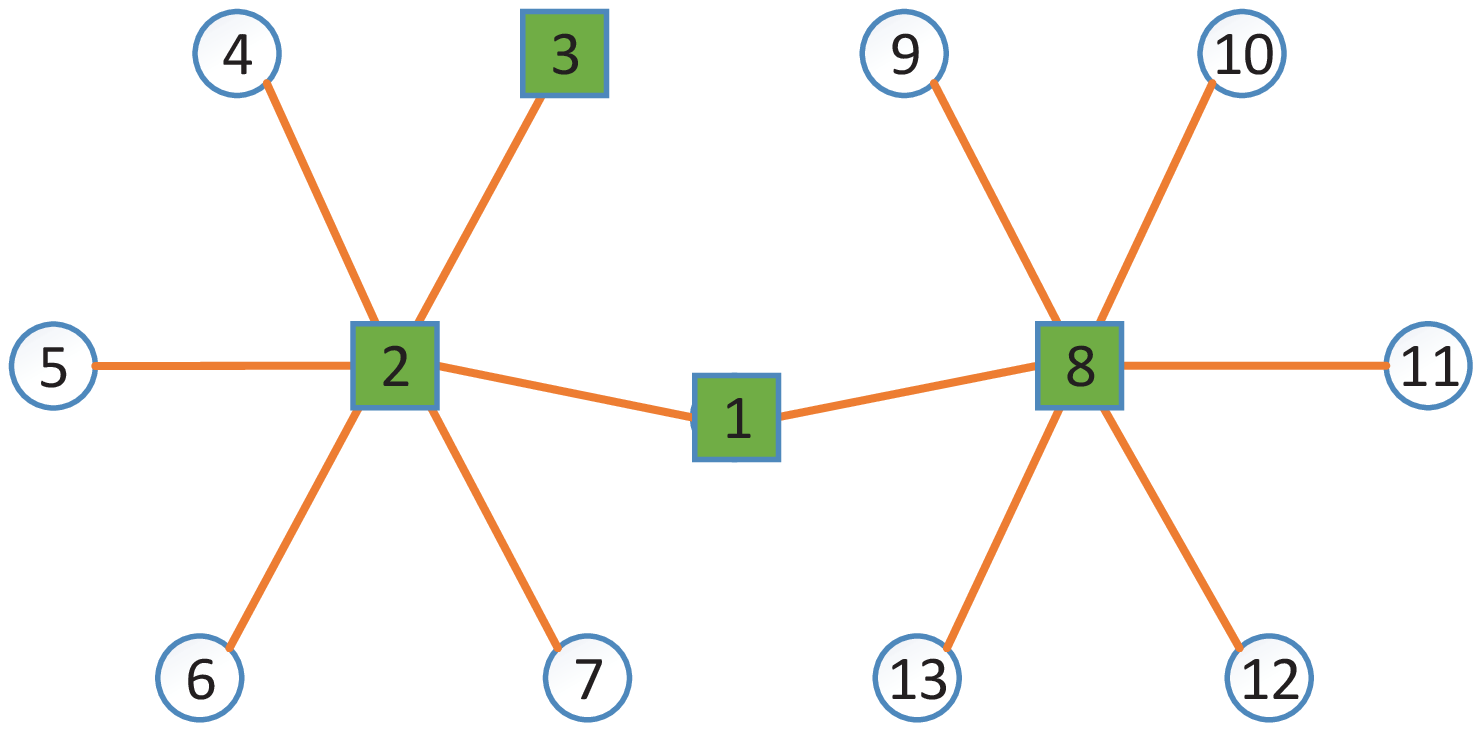}}
  \hspace{0.15in}
  \subfigure[]{
    %\label{fig:subfig:b} %% label for second subfigure
    \includegraphics[width=1.5in]{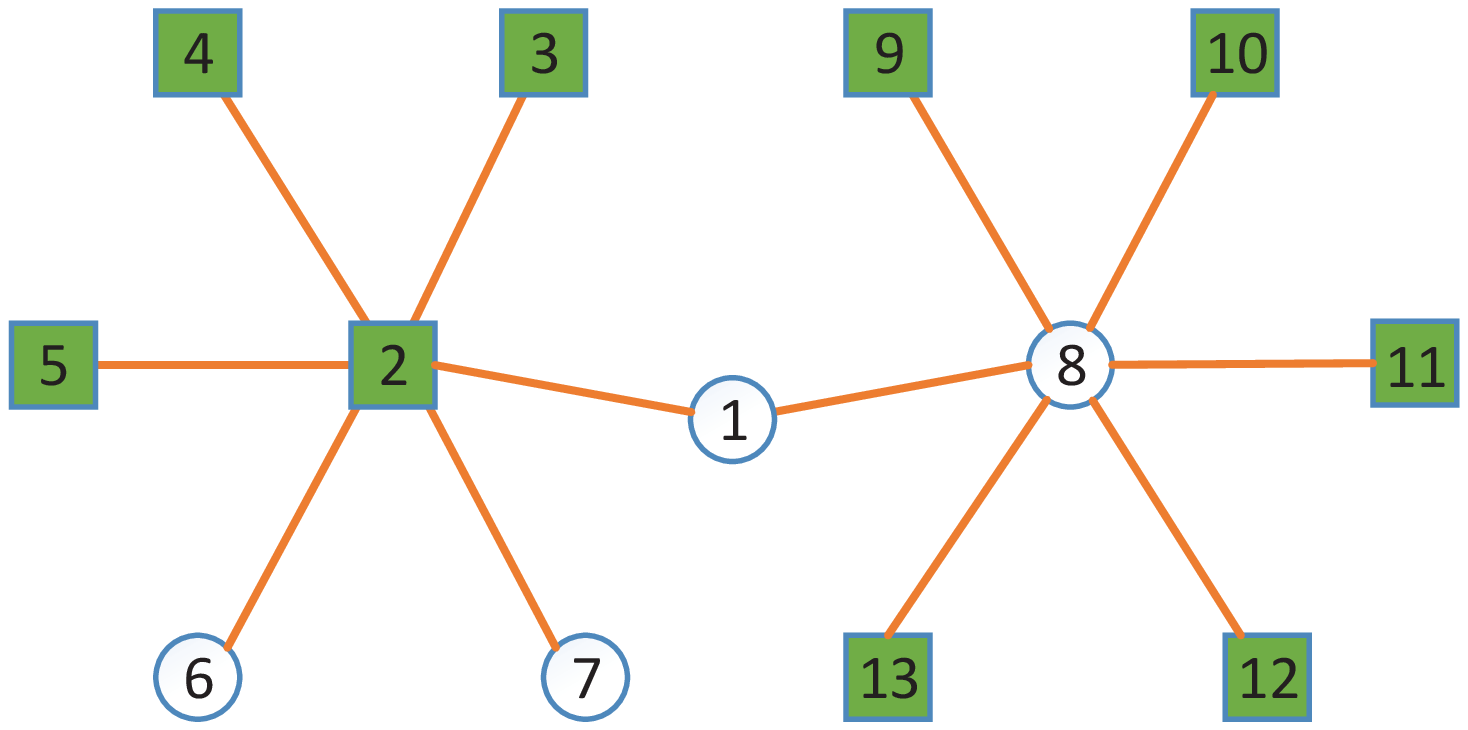}}
  \caption{Several pinning schemes that yield $\max_{|{\mathcal S}|=l}[\lambda_{1}(L(\mathcal S|\mathcal S))]$, when (a) $l=2$; (b) $l=3$; (c) $l=4$; (d) $l=9$.
  A filled square represents a controlled node; a circle represents an uncontrolled node. }
  \label{fig:5} %% label for entire figure
\end{figure}

\end{exmp}

%another example
\begin{exmp}\label{ex:SF_25_nodes}
Consider two networks with 25 nodes as shown in Fig. \ref{fig:9}. The two networks are generated by the Barab\'{a}si-Albert (BA) preferential attachment algorithm \cite{AlBa02}. Starting with three fully-connected nodes, at each time step a new node is added and connected to $m$ existing nodes in the network according to the node-degree preferential attachment scheme \cite{AlBa02}. Figures \ref{fig:9}(a) and (b) show the networks with $m=1$ and $m=2$, respectively.
Figure \ref{fig:7} shows the trends of $\max_{|{\mathcal S}|=l}[\lambda_{1}(L(\mathcal S|\mathcal S))]$ and $\lambda_{l+1}(L_N)$ as $l$ increases. From Fig. \ref{fig:7}, one can see that the Laplacian spectrum of $L_N$ can well estimate $\max_{|{\mathcal S}|=l}[\lambda_{1}(L(\mathcal S|\mathcal S))]$ when $l$ is relatively small ($\frac{l}{N} \leq 50\%$).

\begin{figure}[ht]
\begin{minipage}{0.45\linewidth}
\centerline{\includegraphics[width=1\textwidth]{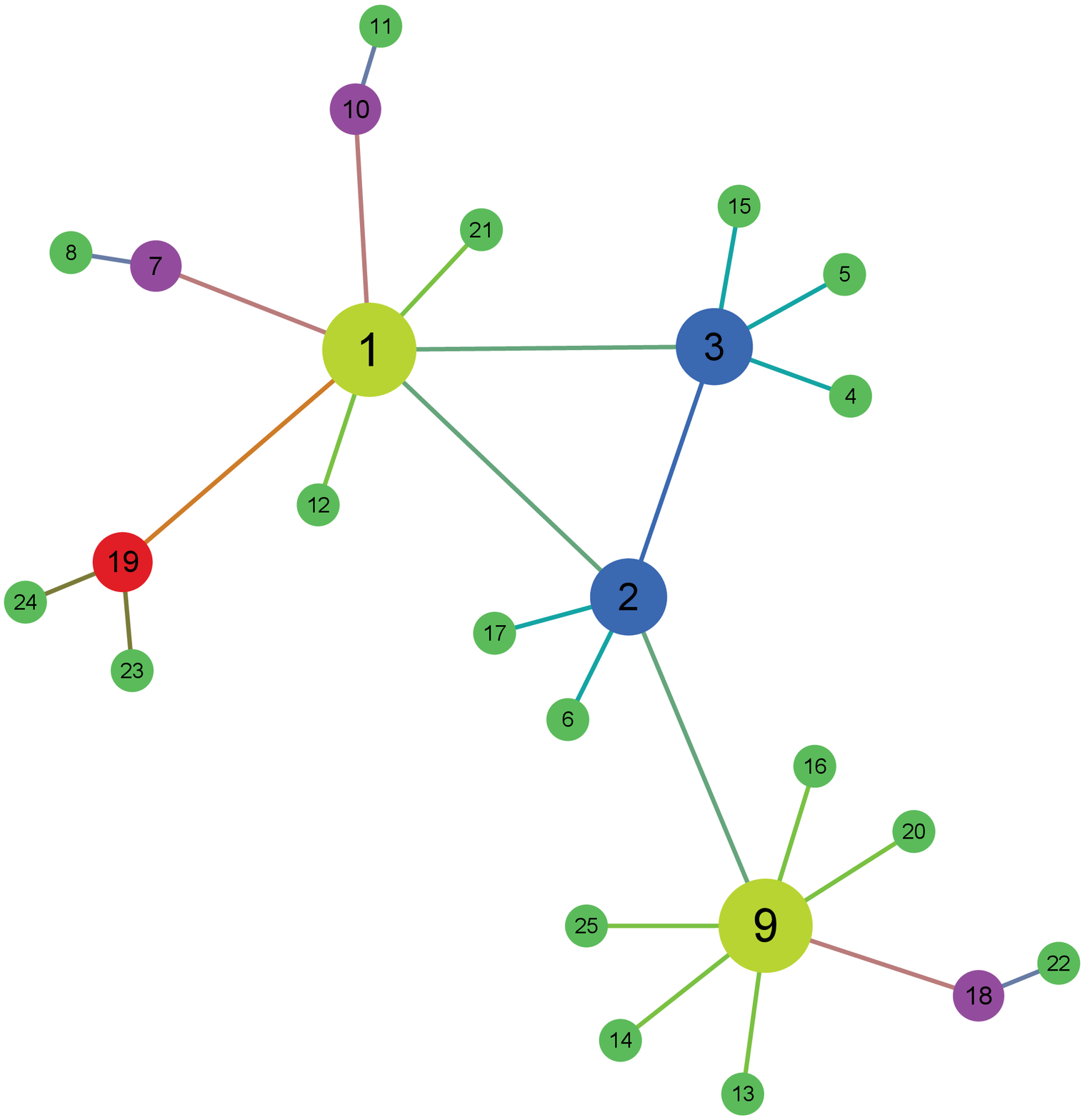}}
\centerline{(a) network I}
\end{minipage}
\qquad
\begin{minipage}{0.45\linewidth}
\centerline{\includegraphics[width=1\textwidth]{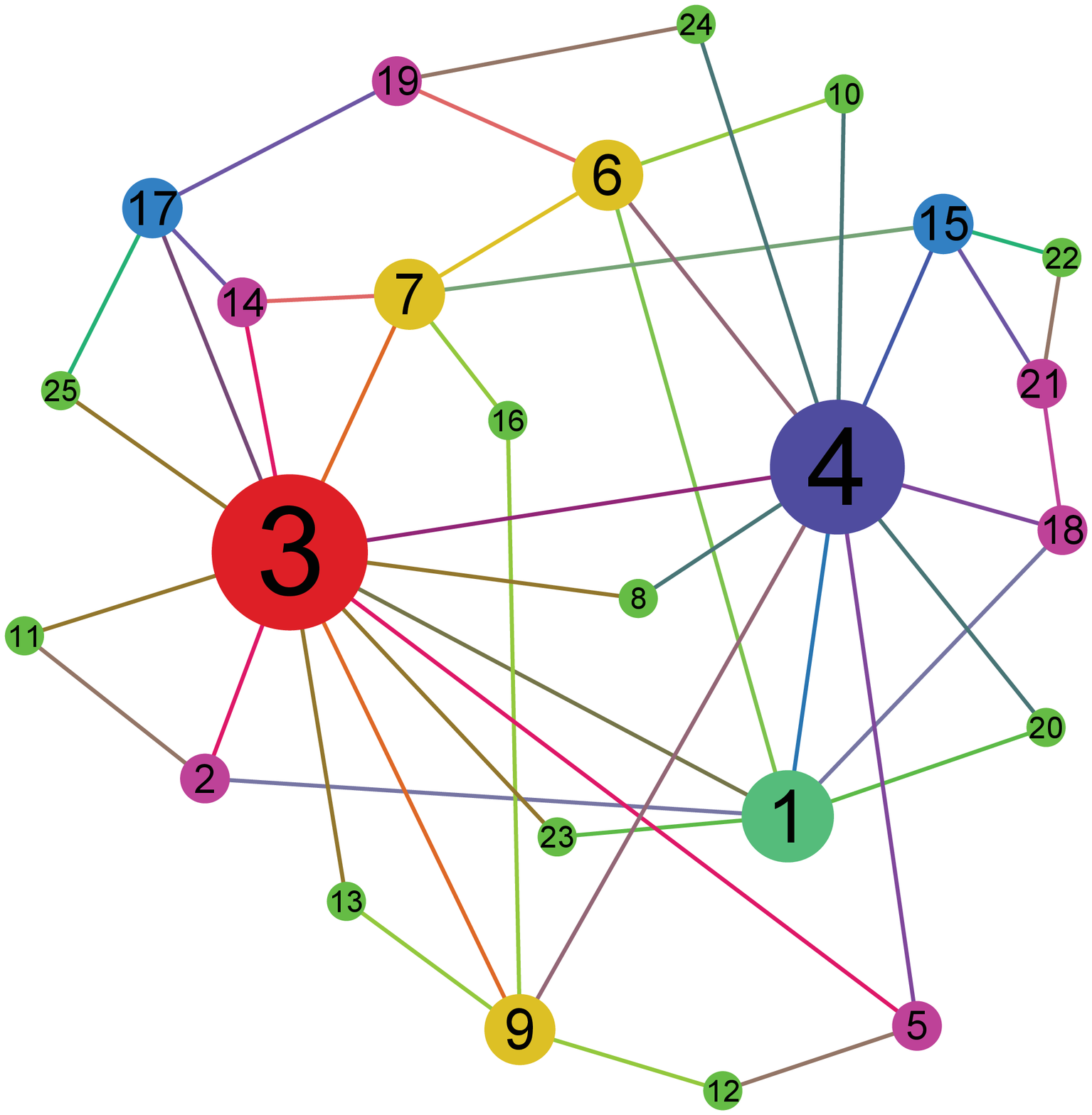}}
\centerline{(b) network II}
\end{minipage}
\caption{Topologies of two networks with 25 nodes each, generated by the BA preferential attachment scheme. (a) $m=1$; (b) $m=2$.}
\label{fig:9}
\end{figure}

\begin{figure}
\centering
\begin{minipage}[t]{0.49\textwidth}
\centering
 \includegraphics[width=3in]{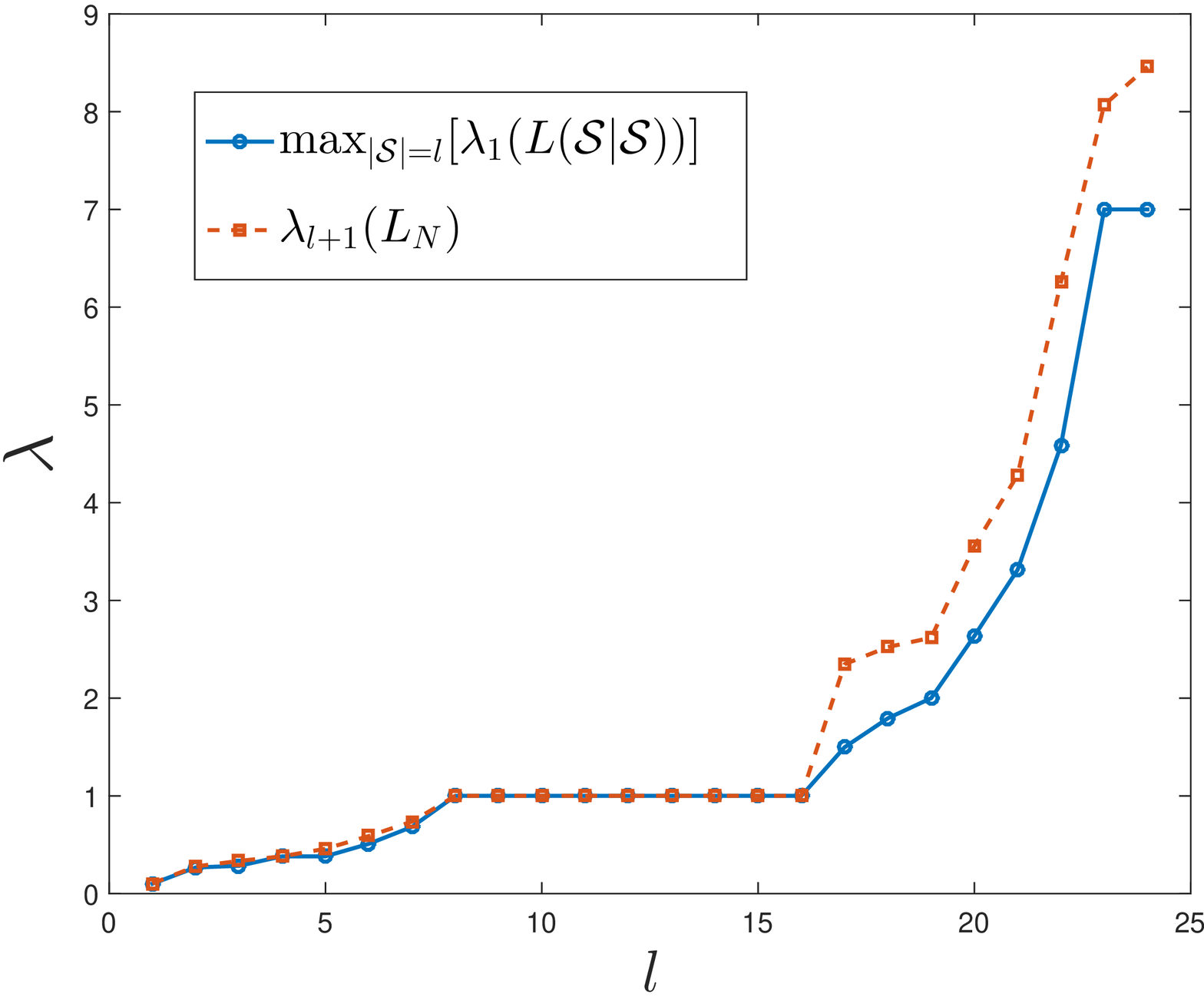}
\centerline{\small (a)}
\end{minipage}
\begin{minipage}[t]{0.49\textwidth}
\centering
 \includegraphics[width=3in]{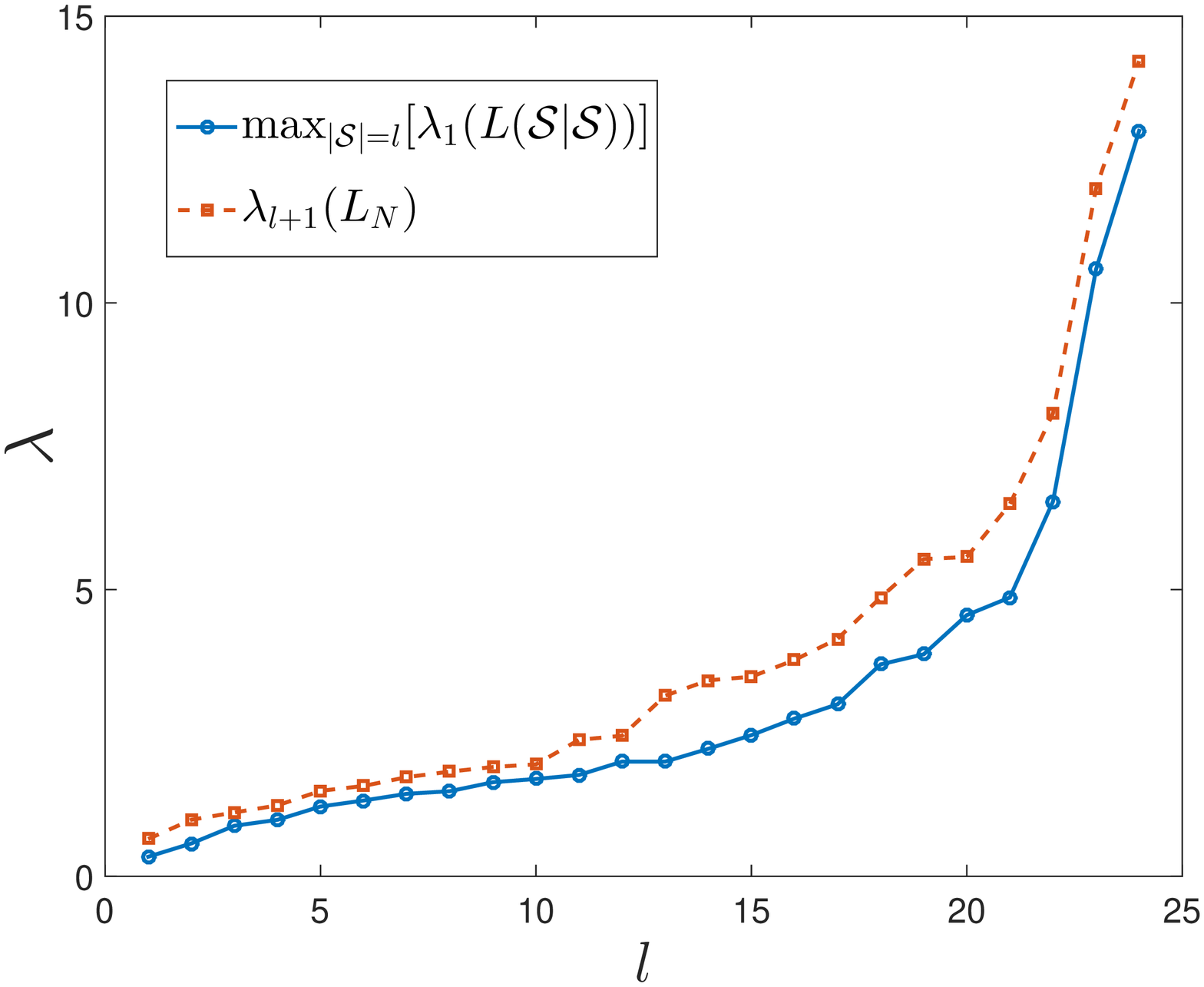}
\centerline{\small (b)}
\end{minipage}
\caption{The trends of $\max_{|{\mathcal S}|=l}[\lambda_{1}(L(\mathcal S|\mathcal S))]$ and $\lambda_{l+1}(L_N)$ as the number $l$ of pinned nodes increases.    (a) for network I in Fig. \ref{fig:9}; (b) for network II in Fig. \ref{fig:9}.}
\label{fig:7}
\end{figure}
\end{exmp}

One can also observe from both Figs. \ref{fig:4} and \ref{fig:7} that $\max_{|{\mathcal S}|=l}[\lambda_{1}(L(\mathcal S|\mathcal S))]$ is non-decreasing as $l$ increases. Actually, this is true and proven by the theorem below.

\begin{theorem}\label{the:mul_eigenvalue_compare}
It holds that
$$\max_{|{\mathcal S}|=l}[\lambda_{1}(L(\mathcal S|\mathcal S))] \geq \max_{|{\mathcal S}|=l-1}[\lambda_{1}(L(\mathcal S|\mathcal S))]$$
for $1\leq l\leq N-1$\,.
\end{theorem}

\begin{IEEEproof}
Suppose that the chosen nodes $s_1, s_2, \dots, s_{l-1}$ can make $\lambda_{1}({L}(\mathcal S_{l-1}|\mathcal S_{l-1}))=\max_{|\mathcal S|=l-1}[\lambda_{1}({L}(\mathcal S|\mathcal S))]$. Besides these $l-1$ nodes, choose one more node $s_l$ to be controlled. Then, the resulting matrix $L(\mathcal S_{l}|\mathcal S_{l})$ satisfies $\lambda_{1}(L(\mathcal S_{l}|\mathcal S_{l})) \geq \lambda_{1}(L(\mathcal S_{l-1}|\mathcal S_{l-1}))$ according to Lemma \ref{lem:eig_com}. Note also that $\max_{|\mathcal S|=l}[\lambda_{1}(L(\mathcal S|\mathcal S))]\geq\lambda_{1}(L(\mathcal S_{l}|\mathcal S_{l}))$. Combining the above three expressions yields that $\max_{|\mathcal S|=l}[\lambda_{1}(L(\mathcal S|\mathcal S))]\geq\max_{|\mathcal S|=l-1}[\lambda_{1}(L(\mathcal S|\mathcal S))]$.
\end{IEEEproof}

%The following corollary can be easily obtained by Theorem \ref{the:mul_eigenvalue_compare}.
\begin{corollary}\label{co:the3}
Let $p,q$ be the numbers of nodes to be controlled, satisfying $1\leq q<p \leq (N-1)$. Then,
\begin{equation*}
\max_{|\mathcal S|=p}[\lambda_{1}(L(\mathcal S|\mathcal S))] \geq \max_{|\mathcal S|=q}[\lambda_{1}(L(\mathcal S|\mathcal S))] \,.
\end{equation*}
\end{corollary}

This corollary implies that $\max_{|\mathcal S|=l}[\lambda_{1}(L(\mathcal S|\mathcal S))]$ is non-descending as $l$ increases.

\begin{theorem}\label{the:mul_upper_k}
It holds that
    \begin{equation}\label{eq:mul_upper_k}
    \lambda_{1}(L(\mathcal S|\mathcal S)) \leq k_{\min}\,,
    \end{equation}
where $k_{\min}$ is the minimal degree of the uncontrolled nodes in the network.
\end{theorem}

\begin{IEEEproof}
Note that the diagonal elements of $L(\mathcal S|\mathcal S)$ are equal to the degrees of the nodes in $\mathcal V\setminus\mathcal S$ correspondingly. Using $\lambda_1(A) \leq \min_{1\leq i\leq N}[a_{ii}]$ in Lemma \ref{lem:eig_diag}, one has that $\lambda_{1}(L(\mathcal S|\mathcal S))\leq \min_{p_j \in \mathcal V\setminus\mathcal S} [k_{p_j}]$.
%, where $p_j$ represents a node in the set $\mathcal V\setminus\mathcal S$ of uncontrolled nodes.
\end{IEEEproof}

\begin{remark}\label{re:mul_upper_k}
When the number of pinned nodes increases, $\lambda_1(L(\mathcal S|\mathcal S))$ usually increases. However, it is bounded by the minimal degree of the uncontrolled nodes, according to Theorem \ref{the:mul_upper_k}. Therefore, at the time that $\lambda_1(L(\mathcal S|\mathcal S))$ reaches $k_{\min}$ in (\ref{eq:mul_upper_k}), it is best to control the node(s) with the smallest degree in the network so as to improve $\lambda_1(L(\mathcal S|\mathcal S))$ if $l$ (i.e.,  $|\mathcal S|$) increases further.
\end{remark}

In the following, the double-star graph in Fig. \ref{fig:3} is used again to verify Theorem \ref{the:mul_upper_k}:\\
{\textit 1)} From Theorem \ref{the:mul_upper_k}, it holds that $\max_{|\mathcal S |=l}[\lambda_1(L(\mathcal S|\mathcal S))] \leq 1$ for $l=1,\dots,9$, because there exists at least one uncontrolled leaf node when $l\leq 9$. This is clear from Fig. \ref{fig:4}.\\
{\textit 2)} Choose $\mathcal V\setminus\mathcal S$ to be $\{2, 8\}$ or $\{2\}$ or $\{8\}$, where $k_2 = k_8 = 6$. There exists at least one uncontrolled node with degree 6. From Theorem \ref{the:mul_upper_k}, it holds that $\lambda_1(L(\mathcal S|\mathcal S))\leq 6$. Actually, this is true because $\lambda_1(L(\mathcal S|\mathcal S)) = 6$.

Finally, another pair of lower and upper bounds for $\lambda_{1}\big(L(\mathcal S|\mathcal S)\big)$ were introduced, which were partially studied in \cite{PiSu16}.

\begin{theorem}\label{the:mul_upper_l}
For $\mathcal S = \{s_1,\dots,s_{l}\}$ and $\mathcal V\setminus\mathcal S = \{p_1,\dots,p_{N-l}\}$, let $w_{p_j}$ be the number of nodes in $\mathcal S$ that are connected with node $p_j \in \mathcal V\setminus\mathcal S$. Then, \\
\begin{equation}\label{eq:12}
\begin{split}
\lambda_{1}\big(L(\mathcal S|\mathcal S)\big) &\leq \frac{w_{p_1}+w_{p_2}+\dots+w_{p_{N-l}}}{N-l} \,;
%&=\frac{k_{s_1}+k_{s_2}+ \cdots +k_{s_{l}}-2e}{N-l}\,,
\end{split}
\end{equation}
\begin{equation}\label{eq:18}
\lambda_1(L(\mathcal S|\mathcal S)) \geq \min_{p_{i}\in \mathcal V\setminus\mathcal S}[w_{p_1},\dots,w_{p_{N-l}}]\,.
\end{equation}
\end{theorem}

\begin{IEEEproof}
\emph{(i)} Here, $\bb H$ is the induced subgraph of $\bb G$ with node set $\mathcal V\setminus\mathcal S$. So, $L(\mathcal S|\mathcal S)=L(\bb H)+\Lambda$, where $\Lambda=\diag\{w_{p_1},\dots,w_{p_{N-l}}\}$. Let $x_0=\frac{1}{\sqrt{N-l}}(1,1,\dots,1)^{\top}$ be an $(N-l)\times 1$ column vector. From Lemma \ref{lem:eigenvalue_definition}, one has
\begin{align}
   \lambda_{1}\big(L(\mathcal S|\mathcal S)\big)
   &=\min_{\|x\|=1}[x^{\top}L(\bb{H})x+x^{\top}\Lambda x]\nonumber\\
   &\leq x_{0}^{\top}L(\bb{H})x_{0}+x_{0}^{\top}\Lambda x_{0}\nonumber\\
   &=\frac{w_{p_1}+w_{p_2}+\dots+w_{p_{N-l}}}{N-l}\,.
   \label{eq:5}
\end{align}

\emph{(ii)} In Lemma \ref{lem:weyl}, take $A=L(\bb H)$, $B=\Lambda$, and $A+B=L(\mathcal S|\mathcal S)$. With $i=j=1$ in (\ref{eq:14}), it follows that
\begin{equation*}
\lambda_1(L(\mathcal S|\mathcal S))\geq \lambda_{1}(L(\bb H))+\lambda_{1}(\Lambda)\,,
\end{equation*}
which is reduced as
$\lambda_1(L(\mathcal S|\mathcal S))\geq \lambda_{1}(\Lambda)$.
Note that $\Lambda=\diag\{w_{p_1},\dots,w_{p_{N-l}}\}$.
Thus, $\lambda_1(L(\mathcal S|\mathcal S)) \geq \min_{p_{i}\in \mathcal V\setminus\mathcal S}[w_{p_1},\dots,w_{p_{N-l}}]$.
The proof is completed.
\end{IEEEproof}

%\begin{remark}\label{re:mul_upper_l-a}
%In $ \frac{k_{s_1}+k_{s_2}+ \cdots +k_{s_{l}}-2e}{N-l}$,  implies that the upper bound for $\lambda_{1}\big(L(\mathcal S|\mathcal S)\big)$ can be improved when one controls the nodes with larger degrees and with relative less connections between the controlled nodes themselves.
%\end{remark}

\begin{remark}
A similar upper bound in (\ref{eq:12}) was also proposed by Theorem 1 of \cite{PiSu16}. Since the notations in \cite{PiSu16} are quite different from those in this paper, for the readers' convenience, a brief and clear proof is provided here.
In fact, a different lower bound that only uses edge information is provided here.
\end{remark}

Consider the situation where only one node,  denoted by $i$, is controlled. From Theorem \ref{the:mul_upper_l}, one obtains the following result.

\begin{corollary}\label{co:one_upper_1}
It holds that
\begin{equation}\label{eq:22}
   \lambda_{1}(L(i|i))\leq\frac{k_i}{N-1} \leq 1 \,,
\end{equation}
where $k_i$ is the degree of the pinned node $i$.
\end{corollary}

Actually, $\lambda_{1}(L(i|i))= 1$ if and only if the node $i$ connects with the rest $N-1$ nodes.
\hspace{0.1cm}\emph{(i)} Necessity: To make $\lambda_{1}(L(i|i))=1$, $k_i$ needs to be at least $N-1$ in (\ref{eq:22}).
\hspace{0.1cm}\emph{(ii)} Sufficiency: Suppose that node $i$ connects with the rest $N-1$ nodes. Then, $\lambda_{1}\big(L(i|i)\big)=\lambda_{1}\big(L(\bb{H})+I_{N-1}\big)=\lambda_{1}(L(\bb{H}))+\lambda_{1}(I_{N-1}) =1$, using Lemma \ref{lem:eigenvalue_polynomial}.
\begin{remark}
Corollary \ref{co:one_upper_1} shows that,
if there exists a node that connects with all the other $N-1$ nodes in the graph, then $\lambda_1(L(i|i))$ reaches its maximum 1 by pinning this node.
\end{remark}

\subsection{The Smallest Eigenvalue $\lambda_1$ of Two Typical Graphs}

The following result gives the smallest eigenvalue $\lambda_{1}\big(L(i|i)\big)$ when pinning one node in the star graph $\mathbb{S}_N$ with $N$ nodes.

\begin{proposition}\label{pro:one_star}
Consider an undirected star graph $\mathbb{S}_N$ with $N>2$ nodes. Suppose that node 1 is the center that has $N-1$ neighbors. Then, \\
{\textit 1)} $\lambda_{1}\big(L(j|j)\big)=\frac{N-\sqrt{N^2-4}}{2}$ for $j=2,\dots,N$;\\
{\textit 2)} $\lambda_{1}\big(L(1|1)\big) = 1$.
% if $\mathcal S = \{1\}$.
\end{proposition}
\begin{IEEEproof}
{\textit 1)} When a leaf node $j$ is pinned, the eigenvalues of $L(j|j)$ can be obtained as follows:
\begin{equation*}
\begin{split}
&|\lambda I_{N-1}-L(j|j)|\\
&=\begin{vmatrix}
     \lambda-(N-1)&1        &1        &\dots&1\\
     1            &\lambda-1&0        &\dots&0\\
     1            &0        &\lambda-1&\dots&0\\
     \dots        &\dots    &\dots    &\dots&\dots\\
     1&0&0&\dots&\lambda-1
\end{vmatrix}_{(N-1)\times(N-1)}\,.\\
\end{split}
\end{equation*}
Multiplying both sides by $(\lambda-1)$ yields
\begin{small}
\begin{equation*}
\begin{split}
&(\lambda-1)|\lambda I_{N-1}-L(j|j)|\\
&=\begin{vmatrix}
\scriptstyle{(\lambda-1)(\lambda-(N-1))}&1        &1        &\dots&1\\
(\lambda-1)             &\lambda-1&0        &\dots&0\\
(\lambda-1)             &0        &\lambda-1&\dots&0\\
     \dots              &\dots    &\dots    &\dots&\dots\\
(\lambda-1)             &0        &0        &\dots&\lambda-1
\end{vmatrix}_{(N-1)\times(N-1)}\\
\end{split}
\end{equation*}
\begin{equation*}
\begin{split}
&=\begin{vmatrix}
\scriptscriptstyle{(\lambda-1)(\lambda-(N-1))-(N-2)}&1                      &1                      &\dots&1\\
0                                                   &\scriptstyle{\lambda-1}&0                      &\dots&0\\
0                                                   &0                      &\scriptstyle{\lambda-1}&\dots&0\\
     \dots                                          &\dots                  &\dots                  &\dots&\dots\\
0                                                   &0                      &0                      &\dots&\scriptstyle{\lambda-1}
\end{vmatrix}_{(N-1)\times(N-1)}\\
&=(\lambda-1)^{N-2}(\lambda^2-N\lambda+1)\,.
\end{split}
\end{equation*}
\end{small}\\
Thus, the eigenvalues of $L(j|j)$ are $\frac{N-\sqrt{N^2-4}}{2},1,\frac{N+\sqrt{N^2-4}}{2}$, so  $\lambda_{1}\big(L(j|j)\big)=\frac{N-\sqrt{N^2-4}}{2}\,.$
\\
{\textit 2)} When the central node 1 is pinned, $L(1|1)$ is the $(N-1)\times(N-1)$ identity matrix. It follows that $\lambda_{1}\big(L(1|1)\big) = 1$. The proof is completed.
\end{IEEEproof}

\begin{remark}
If one leaf node is pinned, $\lambda_{1}\big(L(j|j)\big)$ approximates $1/N$ as the network size $N$ is large, for $j = 2, \dots, N$.
\end{remark}

\begin{proposition}\label{pro:mul_upper_k0}
Let $\mathbb{K}_N$ be a fully-connected network with $N$ nodes, and $\mathcal S = \{s_1,s_2,\dots,s_l\}$ be the set of its controlled nodes, with $1 \leq l \leq N-1$. No matter which nodes are chosen from the node set of $\bb{K}_N$, it holds that
\begin{equation}
\begin{split}
\lambda_1(L(\mathcal S|\mathcal S))=l\,; \quad \lambda_i(L(\mathcal S|\mathcal S))=N,~i = 2,\dots,N-l\,.
\end{split}
\end{equation}
\end{proposition}

\begin{IEEEproof}
{\textit 1)} For $1\leq l<N-1$, one has
\begin{small}
\begin{equation*}
\begin{split}
&|\lambda I_{N-l}-L(\mathcal S|\mathcal S)|\\
&=\begin{vmatrix}
     \lambda-(N-1)&1&\dots&1\\
     1&\lambda-(N-1)&\dots&1\\
     \dots&\dots&\dots&\dots\\
     1&\dots&1&\lambda-(N-1)
\end{vmatrix}_{(N-l)\times(N-l)}\\
\end{split}
\end{equation*}
\begin{equation*}
\begin{split}
&=(\lambda-l)
\begin{vmatrix}
     1&1&\dots&1\\
     1&\lambda-(N-1)&\dots&1\\
     \dots&\dots&\dots&\dots\\
     1&\dots&1&\lambda-(N-1)
\end{vmatrix}_{(N-l)\times(N-l)}\\&=(\lambda-l)
\begin{vmatrix}
     1&1&\dots&1\\
     0&\lambda-N&\dots&0\\
     \dots&\dots&\dots&\dots\\
     0&\dots&0&\lambda-N
\end{vmatrix}_{(N-l)\times(N-l)}\\&=(\lambda-l)(\lambda-N)^{N-l-1}\,.
\end{split}
\end{equation*}
\end{small}
It follows that the eigenvalues of $L(\mathcal S|\mathcal S)$ are $l$ and $N$.
\\
{\textit 2)} For $l=N-1$, it is easy to verify that $L(\mathcal S|\mathcal S)=N-1=l$. This completes the proof.
\end{IEEEproof}

\section{Application of Theoretical Results to Practical Networks}\label{se:app}
%\section{Pinning Control of Practical Networks}\label{se:app}

This section presents applications of the theoretical results obtained in Section \ref{se:main} to  scale-free, small-world, and two real-world networks.

To illustrate the effectiveness of our proposed schemes, they are compared with related works on node-selection schemes. References \cite{SoCa10,WaCh02a,GhAr16} investigated degree-based pinning schemes and preferred to pin the nodes with largest degrees; refs. \cite{RoLiLu09,ZhXi10,JiLiRo11} discussed betweenness centrality (BC)-based pinning schemes.
This work tries to optimize node-selection schemes by analyzing the spectral properties of the grounded Laplacian matrix.
Numerical simulations on typical complex networks in this section show that our proposed schemes have better performance than the aforementioned related works.

\subsection{Scale-Free and Small-World Networks}

Consider scale-free and small-world networks with $N = 1000$.
The scale-free network is generated by the Barab\'{a}si-Albert (BA) preferential attachment algorithm \cite{AlBa02},
which has been explained in Example \ref{ex:SF_25_nodes}.
Start with a fully-connected network with 5 nodes and set $m = 5$. The resulting scale-free network has a power-law degree distribution in theory. The small-world network is generated by the Newman-Watts (NW) algorithm \cite{NeWa99}. Start with a ring, in which each node connects with its $K$ nearest neighbors ($K/2$ on either side). Then, independently add edges with probability $p$ between any pair of unconnected nodes. Here, set $K=4$ and $p=0.006$ in the simulations.

\begin{figure}
\centering
\begin{minipage}[t]{0.49\textwidth}
\centering
 \includegraphics[width=3.4in]{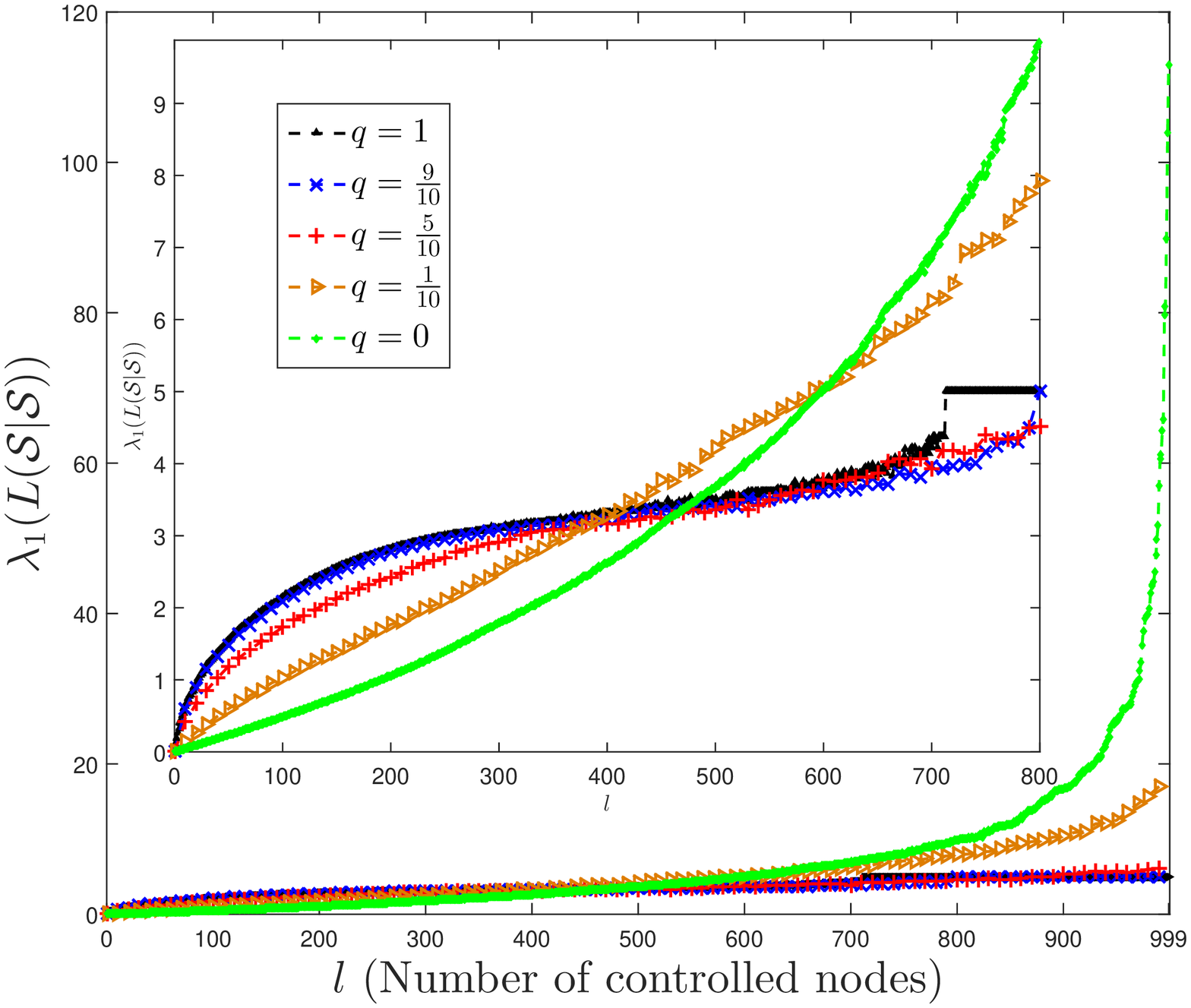}
\centerline{\small (a)}
\end{minipage}
\begin{minipage}[t]{0.49\textwidth}
\centering
 \includegraphics[width=3.4in]{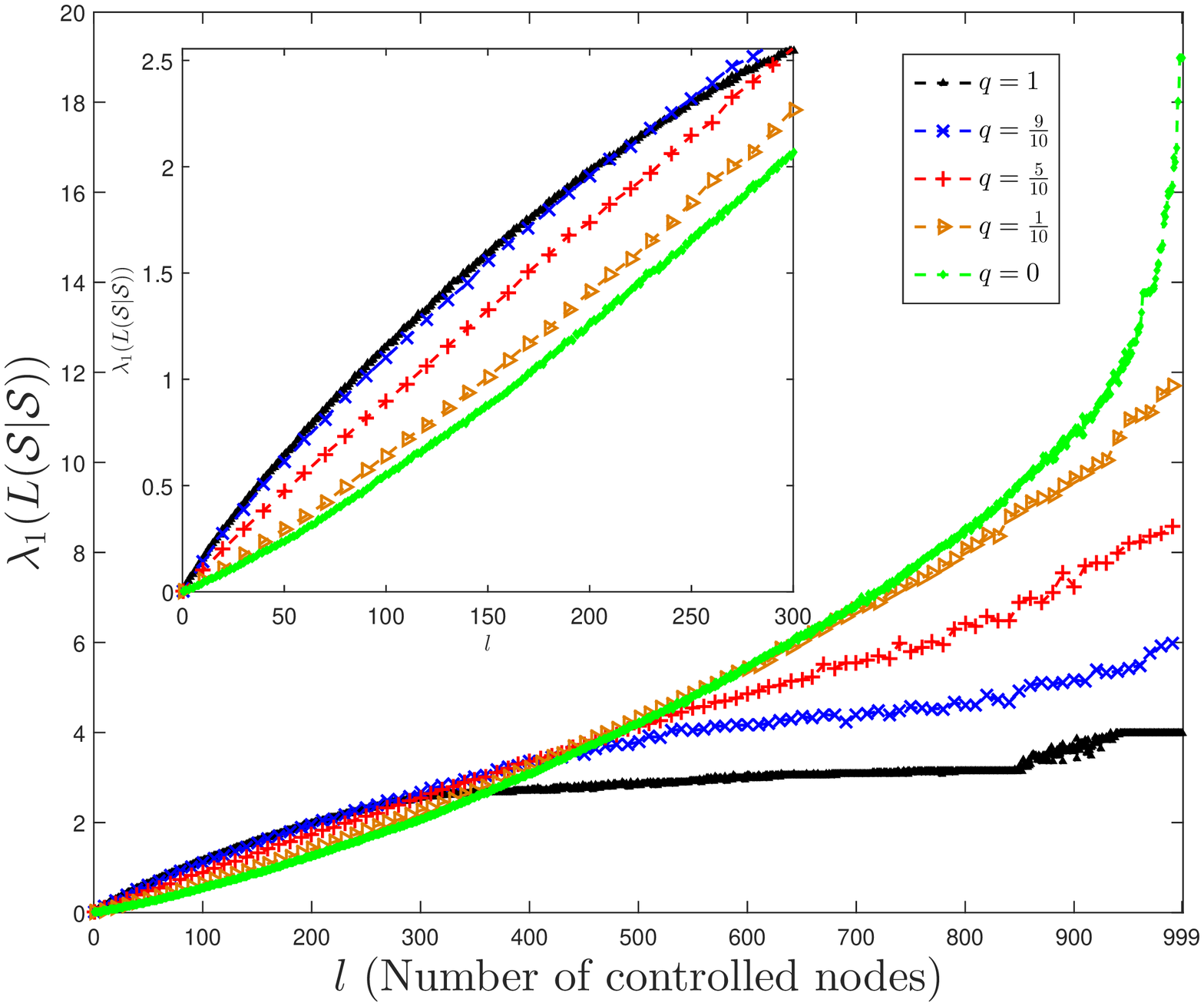}
\centerline{\small (b)}
\end{minipage}
\caption{The value of $\lambda_{1}(L(\mathcal S|\mathcal S))$ versus the number $l$ ($0\leq l\leq 999$) of controlled nodes. The curves with different markers correspond to the trends for different pinning strategies parameterized by $q$. (a) scale-free network; (b) small-world network.}\label{fig:app}
\end{figure}

Numerical results are presented in Fig. \ref{fig:app}, where Fig. \ref{fig:app}(a) shows the trends of $\lambda_{1}(L(\mathcal S|\mathcal S))$ as the number of controlled nodes $l$ increases for the scale-free network, and Fig. \ref{fig:app}(b) shows that for the small-world network. In each subfigure, consider different pinning strategies parameterized by $q\in [0,1]$. Specifically,
select $q\cdot l$ nodes with biggest degrees and the other $(1-q)\cdot l$ nodes with smallest degrees.
In the simulations, $q$ is set to be $1$, $9/10$, $1/2$, $1/10$, and $0$, respectively.
The curves in Fig. \ref{fig:app} with different markers correspond to the trends for the five values of $q$.
Extensive numerical simulations show how $q$ affects $\lambda_1(L(\mathcal S |\mathcal S))$ for both scale-free and small-world networks. For $q = 1$ and 0, simulations are carried out over the interval $l\in [0, 999]$ with step size $\delta l=1$; for $q = 9/10, 5/10$, and $1/10$, over $l\in [0, 990]$ with $\delta l=10$.
Every dot plotted in the subfigures is the average value of five simulation runs by randomly choosing nodes with the same degrees.

From Fig. \ref{fig:app}, one can see that $\lambda_1(L(\mathcal S|\mathcal S))$ increases in general as the number of controlled nodes increases.
For a fixed $l\in [0, 400]$, the larger the proportion $q$, the larger the eigenvalue $\lambda_1(L(\mathcal S|\mathcal S))$. On the contrary, for a fixed $l\in [600, 1000)$, the smaller the $q$, the larger the $\lambda_1(L(\mathcal S|\mathcal S))$.
Note that the degree distribution of the scale-free network should follow a power law, here the maximal degree of nodes is 113 and the minimal degree is 5.
If one prefers to pin the nodes with biggest degrees (i.e., $q = 1$), at least one node with smallest degree will be left uncontrolled. By Theorem \ref{the:mul_upper_k}, one has $\lambda_1(L(\mathcal S|\mathcal S))\leq 5$ if $q = 1$. Similarly, $\lambda_1(L(\mathcal S|\mathcal S))\leq 113$ if $q = 0$. Figure \ref{fig:app}(a) verifies the two inequalities and shows that $\lambda_1(L(\mathcal S|\mathcal S))=113$ for $q=0$ and $\lambda_1(L(\mathcal S|\mathcal S))=5$ for $q=1$, when $|\mathcal S|=999$. For the small-world network, its degree distribution is relatively narrow, where the maximal degree is 19 and the minimal degree is 4. As can be seen from Fig. \ref{fig:app}(b), $\lambda_1(L(\mathcal S|\mathcal S))\leq 4$ if $q = 1$ and $\lambda_1(L(\mathcal S|\mathcal S))\leq 19$ if $q = 0$.

Interestingly, Fig. \ref{fig:app} reveals that for both scale-free and small-world networks, in order to maximize $\lambda_1(L(\mathcal S|\mathcal S))$, it is better to pin the nodes with large degrees when the proportion of pinned nodes is relatively small ($\frac{l}{N} < \frac{400}{1000}$), while it is better to pin nodes with small degrees when the proportion of pinned nodes is large
($\frac{l}{N} > \frac{600}{1000}$). This phenomenon can be explained as follows. Theorem \ref{the:mul_upper_k} shows that $\lambda_1(L(\mathcal S|\mathcal S))$ is upper bounded by $k_{\min}$, where $k_{\min}$ is the minimal degree of uncontrolled nodes. If only a small number of nodes are controlled, $\lambda_1(L(\mathcal S|\mathcal S))$ is normally much smaller than $k_{\min}$. Actually, in this case, those uncontrolled nodes with small degrees do not have much impact on the value of $\lambda_1(L(\mathcal S|\mathcal S))$. As the number of pinned nodes grows, $\lambda_1(L(\mathcal S|\mathcal S))$ also increases but is still bounded by $k_{\min}$.
Therefore, it is reasonable to pin the nodes with small degrees when the proportion of the pinned nodes is large.

\subsection{Dolphin Network}

\begin{figure*}
\centering
\includegraphics[width=4.5in]{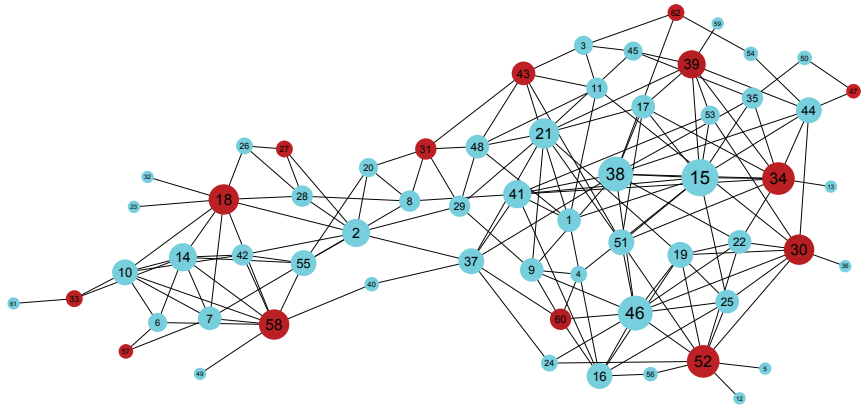}
\caption{Graph of the dolphin network, where the red nodes are selected according to Algorithm \ref{alg:graphpartition} to be controlled.}
\label{fig:dolphins}
\end{figure*}

Consider an undirected and unweighted network of bottlenose dolphins \cite{PiShSu15,Ko17a}, shown in Fig. \ref{fig:dolphins}. The nodes are bottlenose dolphins and an edge indicates a frequent association. There are 62 nodes and 159 edges in Fig. \ref{fig:dolphins}, in which a bigger size of a node implies a larger degree of the node.
Dolphin network has many leaf nodes with degree 1 (14.5\% of the total nodes). According to Theorem \ref{the:mul_upper_k}, $\lambda_1(L(\mathcal S|\mathcal S))$ will be no greater than 1 if there exists a leaf node left uncontrolled. Obviously, controlling all the leaf nodes is not a good pinning strategy, because this does not take into account the different importances of the nodes in a network. According to (\ref{eq:18}) in Theorem \ref{the:mul_upper_l}, it holds that $\lambda_1(L(\mathcal S|\mathcal S)) \geq 1$ if each uncontrolled node has at least one edge connected to a controlled node. Such a graph partition can be realized by a proposed algorithm described as follows:
\begin{algorithm}
  \caption{Find a graph partition such that each node in $\mathcal V \setminus \mathcal S$ has at least one edge connected to set $\mathcal S$. }
  \label{alg:graphpartition}
Input: graph $\mathbb G \Rightarrow {\mathbb G}_0$ and set $\mathcal S = \emptyset $.
\\
\emph{i)} Add all isolated nodes (if exist) to $\mathcal S$. For each connected component in the graph ${\mathbb G}_0$, if there exists node(s) that has links with all the other nodes, randomly choose one of such node(s) and add it to set $\mathcal S$;  if there does not exist such a node in some connected component, for each node with smallest degree, randomly choose one neighbor of the node and add the chosen neighbor node to $\mathcal S$.
\\
\emph{ii)} Remove the newly added controlled nodes and their neighbors from the graph ${\mathbb G}_0$, so as to obtain an induced subgraph with the remaining nodes, which acts as the new graph ${\mathbb G}_1$ for the next iteration;
\\
\emph{iii)} End if the graph is null; otherwise, go to step \emph{i)} with
${\mathbb G}_1 \Rightarrow {\mathbb G}_0$.
\end{algorithm}

Applying Algorithm \ref{alg:graphpartition} to the dolphin network, the controlled node set is obtained as $\mathcal S = \{52,34,18,30,58,39,33,57,27,47,60,62,31,43\}$ with 14 nodes in total. These nodes are marked in red in Fig. \ref{fig:dolphins}. The corresponding $\lambda_1(L(\mathcal S|\mathcal S))=1$.
For comparison, setting the number of controlled nodes to be 14, the values of $\lambda_1(L(\mathcal S|\mathcal S))$ obtained using degree-based and BC-based pinning schemes are shown in Table \ref{tab:pinngscheme_dolphin}.
For the degree-based pinning schemes, choose $^{\rm a}N$ nodes with the largest degrees and $^{\rm b}N$ nodes with the smallest degrees to pin, where ${^{\rm a}N} + {^{\rm b}N}=14$, and then calculate the corresponding values of $\lambda_1(L(\mathcal S|\mathcal S))$.
In the case with multiple choices of node combinations, the value of $\lambda_1(L(\mathcal S|\mathcal S))$ is calculated by taking the average of 100 independent simulation runs.
Numerical simulations are carried out for $^{\rm a}N = 14,7,0$.
For the BC-based pinning scheme,  the $\lambda_1(L(\mathcal S|\mathcal S))$ is calculated when the 14 nodes with largest BC values are chosen to pin.
From Table \ref{tab:pinngscheme_dolphin}, one can observe that these $\lambda_1(L(\mathcal S|\mathcal S))$ under degree-based and BC-based pinning schemes are between $[0.2699, 0.5615]$, much less than 1, which demonstrates that the proposed Algorithm \ref{alg:graphpartition} is more effective.

\begin{table}[!htp]
\centering
\mytoprule
\caption{$\lambda_1(L(\mathcal S|\mathcal S))$ of the dolphin network under different pinning schemes.}\label{tab:pinngscheme_dolphin}
\begin{threeparttable}
\begin{tabular}{    c  |  c  c   c  | c  |c  }
\hline
&\multicolumn{3}{c|}{Degree-based pinning scheme} & \tabincell{c}{BC-based \\pinning \\ scheme } & Algorithm \ref{alg:graphpartition}\\
\hline
$^{\rm a}N$            &14     &7    &0    & $^{\rm BC}|\mathcal S|$	 &\multirow{2}{*}{14} \\
\cline{1-4}
$^{\rm b}N$             &0	    &7	&14  &  $=14$   &               \\	
\hline
$\lambda_1(\mathcal S|\mathcal S)$ &0.5615  &0.5247  &0.2699	 &\tabincell{c}{0.5038}  &1 \\	\hline
\end{tabular}\vspace{8ex}

\mytoprule
\caption{$\lambda_1(L(\mathcal S|\mathcal S))$ of the email network under different pinning schemes.}\label{tab:pinngscheme_email}
\begin{tabular}{    c  |  c  c   c  | c  |c  }
\hline
&\multicolumn{3}{c|}{Degree-based pinning scheme} & \tabincell{c}{BC-based \\pinning \\ scheme } & Algorithm \ref{alg:graphpartition}\\
\hline
$^{\rm a}N$            &266     &133    &0    & $^{\rm BC}|\mathcal S|$	 &\multirow{2}{*}{266} \\
\cline{1-4}
$^{\rm b}N$             &0	    &133	&266  &  $=266$   &               \\	
\hline
$\lambda_1(\mathcal S|\mathcal S)$ &0.3344  &0.3334  &0.3986	 &\tabincell{c}{0.3080}  &1 \\	\hline
\end{tabular}\vspace{1ex}
\footnotesize{$^{\rm a}N$ is the number of controlled nodes selected with the largest degrees.}\\
\footnotesize{$^{\rm b}N$ is the number of controlled nodes selected with the smallest degrees.}\\
\footnotesize{$^{\rm BC}|\mathcal S|$ is the number of controlled nodes selected with the largest betweenness centrality of nodes.}\\
\end{threeparttable}
\end{table}

\subsection{Email Network}

\begin{figure*}
\centering
\includegraphics[width=4.8in]{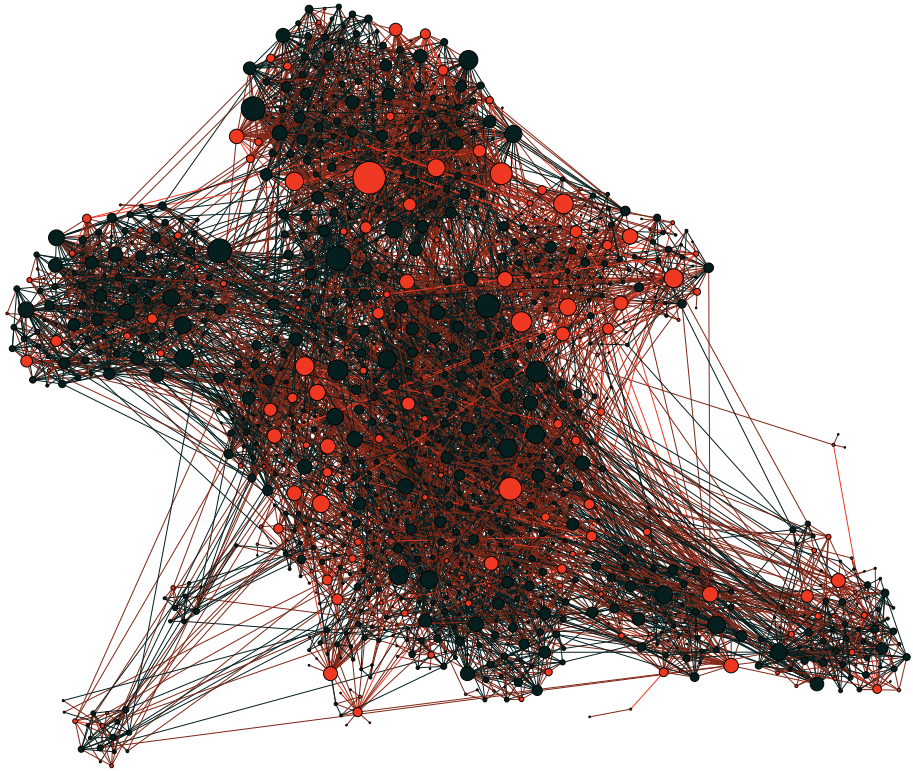}
\caption{Graph of the email network, where the red nodes are selected according to Algorithm \ref{alg:graphpartition} to be controlled.}
\label{fig:email}
\end{figure*}

In this subsection, another real-world network is used to show the effectiveness of Algorithm \ref{alg:graphpartition}.
Fig. \ref{fig:email} shows the email communication network at the University Rovira i Virgili in Tarragona in the south Catalonia in Spain. Nodes are users and each edge represents that at least one email was sent \cite{Ko17b}. The email network consists of 1,133 nodes and 5,453 edges.
In Fig. \ref{fig:email}, there are 151 leaf nodes with degree of 1. According to Algorithm \ref{alg:graphpartition}, a node set contains 266 nodes is obtained to make $\lambda_1(L(\mathcal S|\mathcal S))=1$ by pinning these nodes, which are shown as red nodes in Fig. \ref{fig:email}.
Setting the number of the pinned nodes to be 266, degree-based and BC-based pinning schemes are calculated for comparison with Algorithm \ref{alg:graphpartition}.
The corresponding $\lambda_1(L(\mathcal S|\mathcal S))$ are given in Table II, which shows better performance of Algorithm \ref{alg:graphpartition} than the others.

\section{Conclusions}

An algebraic and graph-theoretic criterion has been introduced for optimizing pinning control of a complex dynamical network. When the node dynamics and the coupling strength of the network are given, the effectiveness of pinning control is measured by the smallest eigenvalue ($\lambda_1 > 0$) of the grounded Laplacian matrix obtained by deleting the rows and columns corresponding to the pinned nodes from the Laplacian matrix of the network. By using tools from graph theory and matrix analysis, formulas have been derived for estimating $\lambda_1$ with a relatively high accuracy using the network topology information. Several important and useful properties of $\lambda_1$ are summarized as follows: \\
%\emph{i)} For an undirected and connected network, $\lambda_1$ is strictly positive provided that at least one node is pinned, which implies that controlling one node can synchronize the entire network.
%\emph{i)} When only one node in a network is pinned, $\lambda_1$ is upper bounded by $1$;  $\lambda_1$ reaches $1$ when the pinned node has connections with all the other nodes in the network. \\
\emph{i)} When $l$ nodes are pinned, $\lambda_1$ is upper bounded by the $(l + 1)$st smallest eigenvalue of the network Laplacian matrix. This property can be used to estimate the minimal number of nodes to be controlled. Especially, when $l$ is relatively small, the Laplacian spectrum of the network can well estimate $\max_{|{\mathcal S}|=l}[\lambda_{1}(L(\mathcal S|\mathcal S))]$.  \\
\emph{ii)} $\lambda_1$ is upper bounded by the minimal degree of the uncontrolled nodes, and also by the average number of edges between the set of  uncontrolled nodes and the set of controlled nodes.

The identified spectral properties have been shown useful and effective for optimizing
pinning control schemes for network synchronization, as demonstrated by illustrative examples.
The designed pinning schemes have be simulated for some typical networks, such as stars, double-star graphs, fully-connected graphs, scale-free and small-world networks.
As future work, it will be important to explore more applications of spectral properties of the grounded Laplacian matrix to network optimization.
Also, it is possible to develop feasible node-selection algorithms by integrating the spectral properties and evolutionary algorithms (e.g. \cite{ZhChChXuHu18,ChZhChTa17}).
%In our ongoing work, we are studying the synthesis of those theoretical estimations in order to develop effective and implementable algorithms for pinning optimization.

\section*{Acknowledgment}                                                             % Place acknowledgements
The authors thank Mr. Wenbin Wei from Wuhan University for helpful discussions.   % here.

\end{document}